\documentclass[10pt]{amsart}
\usepackage{amsmath}
\usepackage{amssymb}
\usepackage{amsfonts}
\usepackage{amsthm,amsxtra}
\usepackage{amsfonts}
\usepackage{nicefrac}

\theoremstyle{plain}
\newtheorem{theorem}{Theorem}[section]
\newtheorem{corollary}{Corollary}[section]
\newtheorem{lemma}{Lemma}[section]
\newtheorem{proposition}{Proposition}[section]
\theoremstyle{definition}
\newtheorem{definition}{Definition}[section]

\theoremstyle{remark}
\newtheorem{remark}{Remark}[section]

\newcommand{\ZZ}{\mathbb Z}
\newcommand{\RR}{\mathbb R}
\newcommand{\CC}{\mathbb C}
\newcommand{\TT}{\mathbb T}

\newcommand{\shalf}{{\scriptstyle \frac{1}{2}}}
\newcommand{\half}{{\frac{1}{2}}}

\begin{document}

\title[Bi-orthogonal Polynomials on the Unit Circle]
{\bf Bi-orthogonal Polynomials on the Unit Circle, regular semi-classical Weights
 and Integrable Systems}

\author{P.J.~Forrester and N.S.~Witte}

\address{Department of Mathematics and Statistics,
University of Melbourne, Victoria 3010, Australia}
\email{\tt p.forrester@ms.unimelb.edu.au; n.witte@ms.unimelb.edu.au}

\begin{abstract}
The theory of bi-orthogonal polynomials on the unit circle is developed for a
general class of weights leading to systems of recurrence relations and derivatives 
of the polynomials and their associated functions, and to functional-difference
equations of certain coefficient functions appearing in the theory.
A natural formulation of the Riemann-Hilbert problem is presented which has
as its solution the above system of bi-orthogonal polynomials and associated functions. 
In particular for the case of regular semi-classical weights on the unit circle
$ w(z) = \prod^m_{j=1}(z-z_j(t))^{\rho_j} $, consisting of $ m \in \mathbb{Z}_{> 0} $
finite singularities, difference equations with respect to the bi-orthogonal polynomial
degree $ n $ (Laguerre-Freud equations or discrete analogs of the Schlesinger equations) 
and differential equations with respect to the 
deformation variables $ z_j(t) $ (Schlesinger equations) are derived completely 
characterising the system. 
\end{abstract}

\maketitle

% \vfill\eject
\section{Introduction}
\setcounter{equation}{0}

The unitary group $U(N)$ with Haar (uniform) measure   
has eigenvalue probability density function (see e.g.~\cite[Chapter 2]{rmt_Fo})
\begin{equation}\label{1.3}
   \frac{1}{(2 \pi )^N N!} \prod_{1 \le j < k \le N} | z_k - z_j |^2,
 \quad z_l := e^{i \theta_l} \in \TT, \quad \theta_l \in (-\pi,\pi] ,
\end{equation} 
where $ \TT = \{z \in \CC: |z|=1 \} $.
Our interest is in averages over $U \in U(N)$ of class functions $ w(U) $ (i.e. symmetric 
functions of the eigenvalues of $ U $ only) which have 
the factorization property $ \prod_{l=1}^N w(z_l) $ for 
$ \{z_1,\dots, z_N \} \in {\rm Spec}(U) $. 
Introducing the Fourier components $\{w_l\}_{l\in \ZZ}$ of the weight
$ w(z) = \sum_{l=-\infty}^\infty w_l z^l $, due to
the well known identity
\cite{ops_Sz}
\begin{equation}\label{1.3a}
  \Big \langle \prod_{l=1}^N w(z_l) \Big \rangle_{U(N)} =
  \det[ w_{i-j} ]_{i,j=1,\dots,N}, 
\end{equation}
we are equivalently studying Toeplitz determinants. 
However we are interested in the situation where the weights are not necessarily positive or 
even real 
valued $ \overline{w(z)} \neq w(z) $ for $ z \in \TT $, where the bar denotes the complex 
conjugate,
and consequently the Toeplitz matrices are non-Hermitian $ \bar{w}_n \neq w_{-n} $. 
The motivations for studying these types of weights comes from many diverse applications,
in particular the gap probabilities and characteristic polynomial averages in random matrix
theory, the spin-spin correlations of the square lattice Ising model, the density 
matrix of one-dimensional systems of impenetrable bosons and
probability distributions for various classes of non-intersecting lattice path problems.
All the above applications correspond to special cases of the weight
\begin{equation}
   w(z) = t^{-\mu}z^{-\omega-\mu}(1+z)^{2\omega_1}(1+tz)^{2\mu} \;
   \begin{cases}
     1, \quad \theta \in (-\pi,\pi-\phi) \cr
     1-\xi, \quad \theta \in (\pi-\phi,\pi]
   \end{cases} ,
\label{1.3b}
\end{equation}
where $ \mu $, $ \omega=\omega_1+i\omega_2 $ are complex parameters and $ \xi $, 
$ t=e^{i\phi} $ are 
complex variables. We have previously shown \cite{FW_2002b} that (\ref{1.3a}) with 
weight (\ref{1.3b}) can, as a function of $ t $, be characterised as a $ \tau$-function 
for the Painlev\'e VI system, 
and thus be expressed in terms of a solution of the sixth Painlev\'e equation.
Here we will develop theory which allows for a characterisation of a general class of
Toeplitz determinants as a function of the discrete variable $ N $. In the particular case 
of the weight (\ref{1.3b}) the resulting difference equations in $ N $ will be studied in
a companion paper \cite{FW_2004b} where it is shown they are equivalent to the 
discrete fifth Painlev\'e equation {${\rm dP}_{\rm V}$}
associated with the degeneration of the rational surface $ D^{(1)}_4 \to D^{(1)}_5 $.

To characterise the unitary group averages (\ref{1.3a})
we have found it necessary to substantially develop the theory of
bi-orthogonal polynomial systems on the unit circle 
$\{ \phi_n(z), \bar{\phi}_n(z) \}^{\infty}_{n=0} $, generalising the Szeg\"o and Geronimus 
theory applying to orthogonal polynomial systems. They are defined by
\begin{equation}
  \int_{\TT} \frac{d\zeta}{2\pi i\zeta} w(\zeta)\phi_m(\zeta)\bar{\phi}_n(\bar{\zeta})  
   = \delta_{m,n} .
\end{equation}
Issues relating to the existence of such systems and their basic properties is 
taken up in Section 2.
A study of bi-orthogonal polynomial systems on the unit circle
was begun by Baxter \cite{Ba_1960,Ba_1961} where their elementary properties were 
elucidated and the connection between the absolute convergence of the trigonometric 
expansion of $ \log w(e^{i\theta}) $ and absolute convergence of certain coefficients
of the bi-orthogonal polynomials was investigated.
Sometime later other useful elements of the theory of orthogonal polynomial systems 
were added by Jones, Nj{\aa}stad and Thron in their study \cite{JNT_1989} of the trigonometric 
moment problem and related quadrature formulae on the unit circle, namely the
associated polynomials and the Hilbert transform of the weight.
The task of determining the functional, difference and differential equations for
orthogonal polynomial systems on the unit circle was initiated in \cite{IW_2001} but was 
not completed nor extended to cover the bi-orthogonal situation.
In extending the theory in this direction we have sought to cast it in a way so that
the appearance is as close to the Szeg\"o-Geronimus theory as possible, and indeed we
find that virtually all formulae from the older theory can be taken over but now
the complex conjugated variables have to be re-interpreted as independent variables. 
We wish to emphasise that the works described above, which are directly relevant to the
present study of certain classes of weights, constitute only a small fraction of existing 
literature on orthogonal polynomial systems on the unit circle. The themes taken up in 
this broader body of work, for example as recounted in \cite{Si_2004}, are really rather 
different problems and are formulated for the most general types of real, positive measures 
on the unit circle. Finally we note that
to a large extent the task of characterising the analogous Hankel determinants has been 
completed for weights defined on the real line and this could be founded upon the standard
theory of orthogonal polynomial systems. In particular we note the works of Bauldry \cite{Ba_1990}, 
Bonan and Clark \cite{BC_1990}, Belmehdi and Ronveaux \cite{BR_1994}, and
Magnus \cite{Ma_1999,Ma_1995a,Ma_1994} where in the last three works the differential and 
difference structures have been clearly revealed and linked to isomonodromy preserving
deformations, culminating in the identification of the simplest nontrivial case with the 
sixth Painlev\'e system. 

In our formulation of the problem we start with a weight, possessing certain 
features to be made precise later, which we take as given and wish to calculate 
the Toeplitz determinants via
auxiliary quantities arising in the bi-orthogonal polynomial theory.
To this end we have derived closed systems of differential relations for the 
polynomials $\{ \phi_n(z) \}^{\infty}_{n=0} $, their reciprocal polynomials 
$\{ \phi^*_n(z) \}^{\infty}_{n=0} $, and associated functions 
$\{ \epsilon_n(z) \}^{\infty}_{n=0} $, $\{ \epsilon^*_n(z) \}^{\infty}_{n=0} $
in Proposition \ref{ops_spectralD}.
In the notation of (\ref{ops_Ydefn}) and Corollary 2.3 let
\begin{equation*}
   Y_n(z;t) := 
   \begin{pmatrix}
          \phi_n(z)   &  \epsilon_n(z)/w(z) \cr
          \phi^*_n(z) & -\epsilon^*_n(z)/w(z) \cr
   \end{pmatrix} .
\end{equation*}
We show
\begin{equation*}
   \frac{d}{dz}Y_{n} := A_n Y_{n} ,
\end{equation*} 
where the entries in the matrix $ A_n $ are parameterised  by four coefficient functions
$ \Omega_n(z) $, $ \Omega^*_n(z) $ , $ \Theta_n(z) $ , $ \Theta^*_n(z) $ in (\ref{ops_YzDer}),
and complete sets of difference and functional relations for these coefficient 
functions are given in Proposition \ref{ops_Linear1} and Corollary \ref{ops_Linear2}. 
We also formulate a $ 2 \times 2$ matrix Riemann-Hilbert problem in Proposition
\ref{ops_RHP} for general classes of weights which parallels the case for orthogonal 
polynomials on the line \cite{IKF_1991,FIK_1991,FIK_1992,De_1999} and whose 
solution is simply related to $ Y_n $. 

For our particular applications the weights are members of the regular 
semi-classical class 
\begin{equation*}
  w(z) = \prod^m_{j=1}(z-z_j(t))^{\rho_j}, \quad \rho_j \in \CC,
\end{equation*}
with an arbitrary number $ m $ of isolated finite singularities located at $ z_j(t) $. 
These are also known as generalised Jacobi weights. Using the projection of the unit
circle onto the interval $ [-1,1] $ a system of orthogonal polynomials on the unit
circle is equivalent to two related systems of orthogonal polynomials on this interval.
The generalised Jacobi orthogonal polynomial systems with support $ [-1,1] $ have been
studied from a number of points of view. For example the Stieltjes electrostatic problem 
was generalised to include an arbitrary number of fixed charges in
\cite{Ge_1958}. Uniform asymptotics of the orthogonal polynomials were first derived 
in \cite{Ba_1983a}, \cite{Ba_1983b} and this was extended to their derivatives in
\cite{Ba_1992}, \cite{Ve_1999a}, \cite{Ve_1999b} and \cite{Ve_2001}. Further properties
of the polynomials, such as their upper bounds, were investigated in \cite{EMN_1994},
\cite{NEM_1994} and estimates of the associated functions and polynomials were made
in \cite{CMN_1991}.
Interest in issues relating to quadrature problems were taken up in \cite{EN_1992}
where lower and upper bounds of the Christoffel function were found and quadrature
inequalities were derived in \cite{MV_1997}. The zeros of the orthogonal polynomials have
also been studied, in particular the spacing of consecutive zeros in \cite{EN_1992} and
their asymptotic formulae in \cite{Ve_1997}. We wish to mention too a work
similar in spirit to our own \cite{GN_1982}, where the orthogonal polynomial system 
$ \{p_n(t)\}^{\infty}_{n=0} $ was defined by
\begin{equation*}
   \int_{C} dt\, \omega(t)p_n(t)t^k = 0,  \quad k=0,\ldots,n-1 .
\end{equation*}
In this work the weight had $ m=3 $ distinct finite singularities of the regular semi-classical
class and its support was a closed curve $ C $
enclosing all the singularities. Recurrence relations in $ n={\rm deg}(p_n) $ for the three-term 
recurrence coefficients were found and heuristic arguments were given for the asymptotic
behaviour of these as well as the polynomials themselves.

A key feature of regular semi-classical weights is that
\begin{equation*}
  \frac{1}{w(z)}\frac{d}{dz}w(z) = \frac{2V(z)}{W(z)} ,
\end{equation*}
where $ V $ and $ W $ are polynomials such that $ {\rm deg}V(z) < m, {\rm deg}W(z)=m $.
The coefficient functions for regular semi-classical weights are polynomials of 
$ z $ with bounded degree, 
$ {\rm deg}\Omega_n(z)={\rm deg}\Omega^*_n(z)=m-1 $, 
$ {\rm deg}\Theta_n(z)={\rm deg}\Theta^*_n(z)=m-2 $
(see Proposition \ref{ops_SCpoly}).
In addition evaluations of these functions at the singular points satisfy bilinear 
relations (see Proposition \ref{ops_Bilinear}) which lead directly to one of the 
pair of coupled discrete Painlev\'e equations. 
Deformation derivatives of the linear system of differential equations above 
with respect to arbitrary trajectories of the finite singularities are given in Proposition
\ref{ops_deformD} which can summarised as
\begin{equation*}
   \frac{d}{dt}Y_{n} := B_n Y_{n}
   = \left\{ B_{\infty} - \sum^{m}_{j=1}\frac{A_{nj}}{z-z_j}\frac{d}{dt}z_j
        \right\} Y_{n} ,\quad \text{where} \quad
   A_n = \sum^{m}_{j=1}\frac{A_{nj}}{z-z_j} ,
\end{equation*} 
and consequently systems of Schlesinger equations for the elements of $ A_{nj} $
(or the coefficient functions evaluated at $ z_j $) are given in 
(\ref{ops_Schl:a}-\ref{ops_Schl:c}). 
It is quite natural that systems governed by regular
semi-classical weights preserve the monodromy data of the solutions $ Y_n $ about 
each singularity $ z_j $ with respect to arbitrary deformations.

We mention at this point that the full definition of a regular semi-classical weight 
given below in Definition \ref{rSC_defn} is restrictive, and has been generalised by 
relaxing some of 
the conditions in a series of works \cite{Ma_1987}, \cite{MR_1992}, \cite{MR_1998}.
In these works the orthogonal polynomial systems were characterised by integral 
representations of semi-classical linear functionals with respect to certain paths 
in the complex plane. An irregular semi-classical weight arising under less
restrictive conditions can be recovered from a particular regular semi-classical weight
through limiting processes involving the coalescence of singular points in the same way as
the fifth to the first Painlev\'e systems are recovered from the sixth.
The differential and difference structures of
orthogonal and bi-orthogonal polynomial systems with polynomial log-derivatives
of the weight ($ {\rm deg}(W)=0 $, $ {\rm deg}(V)>2 $) defined on $ \RR $ have been studied 
\cite{BEH_2003a}, \cite{BEH_2003b} in the context of matrix models. Analogous
structures for a bi-orthogonal polynomial system on the unit circle with a simple
Laurent polynomial log-derivative for the weight were found in \cite{Hi_1996}, motivated
by applications to unitary matrix models.

In Section 2 we derive systems of differential-difference
and functional relations for bi-orthogonal polynomials and associated functions
on the unit circle for a general class of weights and formulate the Riemann-Hilbert 
problem. 
In Section 3 we specialise to regular semi-classical weights and derive bilinear
difference equations. In addition we calculate the deformation derivatives of the 
bi-orthogonal polynomial system, derive a system of Schlesinger equations
and show the deformations are of the isomonodromic type. 

\section{Bi-orthogonal Polynomials on the Unit Circle and Riemann-Hilbert Problem}
\label{OPSsection}
\setcounter{equation}{0}

We consider a complex function for our weight $ w(z) $, analytic in the cut
complex $ z $-plane and which possesses a Fourier expansion
\begin{equation}
  w(z) = \sum_{k=-\infty}^{\infty} w_{k}z^k, \quad
  w_{k} = \int_{\TT} \frac{d\zeta}{2\pi i\zeta} w(\zeta)\zeta^{-k},
\end{equation}
where $ z \in D \subset \CC $ and $ \TT $ denotes the unit circle $ |\zeta|=1 $ with 
$ \zeta=e^{i\theta}, \theta \in (-\pi,\pi] $. 
Hereafter we will assume that $ z^jw(z), z^jw'(z) \in L^1(\TT) $ for all 
$ j \in \ZZ $. We will also assume that the trigonometric sum converges in an annulus
$ D=\{z\in \CC:\Delta_1 < |z| < \Delta_2 \} $ and $ \TT \subset D $. Until we arrive at
our specific weights of interest in Section 3, namely the regular semi-classical class,
we will formally assume convergence holds.  
The doubly infinite 
sequence $ \{ w_k \}^{\infty}_{k=-\infty} $ are the trigonometric moments of the 
distribution $ w(e^{i\theta})d\theta/2\pi $ and define the trigonometric moment 
problem. Define the Toeplitz determinants
\begin{align}
   I^{\epsilon}_{n}[w] 
  & := \det \left[ \int_{\TT} \frac{d\zeta}{2\pi i\zeta} w(\zeta)\zeta^{\epsilon-j+k}
           \right]_{0 \leq j,k \leq n-1},
  \nonumber \\
  & = \det \left[ w_{-\epsilon+j-k} \right]_{0 \leq j,k \leq n-1},
  \nonumber \\
  & = \frac{1}{n!} \int_{\TT^n}
      \prod^{n}_{l=1} \frac{d\zeta_l}{2\pi i\zeta_l} w(\zeta_l)\zeta_l^{\epsilon}
      \prod_{1 \leq j<k \leq n} |\zeta_{j}-\zeta_{k}|^2,
\end{align}
where $ \epsilon $ will take the integer values $ 0,\pm 1 $.

We define a system of bi-orthogonal polynomials 
$ \{ \phi_n(z),\bar{\phi}_n(z) \}^{\infty}_{n=0} $ with respect to the 
weight $ w(z) $ on the unit circle by the orthogonality relation
\begin{equation}
  \int_{\TT} \frac{d\zeta}{2\pi i\zeta} w(\zeta)\phi_m(\zeta)\bar{\phi}_n(\bar{\zeta})  
   = \delta_{m,n} .
\label{ops_onorm}
\end{equation}
This system is taken to be orthonormal and the coefficients in a monomial basis 
are defined by
\begin{align}
   \phi_n(z)
   & = \kappa_n z^n + l_n z^{n-1}+ m_n z^{n-2} + \ldots + \phi_n(0)
     = \sum^{n}_{j=0} c_{n,j}z^j,
   \\
   \bar{\phi}_n(z)
   & = \bar{\kappa}_n z^n + \bar{l}_n z^{n-1}+ \bar{m}_n z^{n-2} + \ldots + \bar{\phi}_n(0)
     = \sum^{n}_{j=0} \bar{c}_{n,j}z^j,
\end{align}
where $ \bar{\kappa}_n $ is chosen to be equal to $ \kappa_n $ without loss of generality
(this has the effect of rendering many results formally identical to the pre-existing
theory of orthogonal polynomials). Notwithstanding the notation $ \bar{c}_{n,j} $ 
in general is not equal to the complex conjugate of $ c_{n,j} $ and is independent of it.
We also define the reciprocal polynomial by
\begin{equation}
   \phi^{*}_n(z) := z^n\bar{\phi}_n(1/z) = \sum^{n}_{j=0} \bar{c}_{n,j}z^{n-j} .
\end{equation}
The bi-orthogonal polynomials can be defined up to an overall factor by the orthogonality
with respect to the monomials
\begin{equation}
  \int_{\TT} \frac{d\zeta}{2\pi i\zeta} w(\zeta)\phi_n(\zeta)\overline{\zeta^j}  
   = 0 \qquad 0 \leq j \leq n-1 ,
\label{ops_orthog:a}
\end{equation}
whereas their reciprocal polynomials can be similarly defined by
\begin{equation}
  \int_{\TT} \frac{d\zeta}{2\pi i\zeta} w(\zeta)\phi^*_n(\zeta)\overline{\zeta^j}  
   = 0 \qquad 1 \leq j \leq n .
\label{ops_orthog:b}
\end{equation}

The linear system of equations for the coefficients $ c_{n,j}, \bar{c}_{n,j} $
arising from 
\begin{align}
  \bar{c}_{n,n}\int_{\TT} \frac{d\zeta}{2\pi i\zeta} w(\zeta)\phi_n(\zeta)\bar{\zeta}^m  
 & = \begin{cases} 0 \quad 0 \leq m \leq n-1 \\ 1 \quad m=n \end{cases} , 
  \label{ops_CoeffEqn:a}\\
  c_{n,n}\int_{\TT} \frac{d\zeta}{2\pi i\zeta} w(\zeta)\zeta^m\bar{\phi}_n(\bar{\zeta})
 & = \begin{cases} 0 \quad 0 \leq m \leq n-1 \\ 1 \quad m=n \end{cases} ,
  \label{ops_CoeffEqn:b}
\end{align}
has the solution
\begin{align}
  c_{nj} & = \frac{1}{\bar{c}_{n,n}}
         \frac{\det \begin{pmatrix}
                     w_{0} 	& \ldots & 0 	& \ldots	& w_{-n} \cr
                     \vdots	& \vdots & \vdots & \vdots 	& \vdots \cr
                     w_{n-1}	& \ldots & 0 	& \ldots 	& w_{-1} \cr
                     w_{n} 	& \ldots & 1 	& \ldots 	& w_{0}  \cr
                    \end{pmatrix} }
              {\det \begin{pmatrix}
                     w_{0} & \ldots & w_{-n} \cr
                     \vdots& \vdots & \vdots \cr
                     w_{n} & \ldots & w_{0}  \cr
                    \end{pmatrix} } ,
  \\
  \bar{c}_{nj} & = \frac{1}{c_{n,n}}
         \frac{\det \begin{pmatrix}
                     w_{0} 	& \ldots & 0 	& \ldots 	& w_{n}  \cr
                     \vdots	& \vdots & \vdots & \vdots 	& \vdots \cr
                     w_{-n+1} 	& \ldots & 0 	& \ldots 	& w_{1}  \cr
                     w_{-n} 	& \ldots & 1 	& \ldots 	& w_{0}  \cr
                    \end{pmatrix} }
              {\det \begin{pmatrix}
                     w_{0} & \ldots & w_{n} \cr
                     \vdots& \vdots & \vdots   \cr
                     w_{-n}& \ldots & w_{0} \cr
                    \end{pmatrix} } ,
\end{align}
and in particular one has the following results.
\begin{proposition}[\cite{Ba_1960}]
The leading and trailing coefficients of the polynomials $ \phi_n(z) $, $ \bar{\phi}_n(z) $
are 
\begin{gather}
   c_{nn} =  \bar{c}_{nn} = \kappa_{n} 
          = \frac{1}{\kappa_n}\frac{I^{0}_{n}}{I^{0}_{n+1}}, \\
   c_{n0} = \phi_{n}(0) = (-1)^n\frac{1}{\kappa_n}\frac{I^{1}_{n}}{I^{0}_{n+1}}, \quad
   \bar{c}_{n0} = \bar{\phi}_{n}(0) 
          = (-1)^n\frac{1}{\kappa_n}\frac{I^{-1}_{n}}{I^{0}_{n+1}} .
\end{gather}
\end{proposition}
\noindent

\begin{proposition}[\cite{Ba_1960}]
The bi-orthogonal system $ \{\phi_{n}(z),\phi^*_{n}(z)\}^{\infty}_{n=0} $ exists
if and only if $ I^{0}_{n} \neq 0 $ for all $ n \geq 1 $. 
\end{proposition}
\begin{remark}
We shall see that a failure of this condition can and generically must occur in the 
case of regular semi-classical weights which contain deformation parameters 
$ z_1, \ldots, z_m $ and this is precisely the condition for
a movable singularity (in this case a pole) in the dynamics with respect to the 
$ z_j $.
\end{remark}
A consequence of the above solutions are the following determinantal and integral 
representations for the polynomials,
\begin{align}
  \phi_{n}(z) & = \frac{\kappa_n}{I^{0}_{n}}
         \det \begin{pmatrix}
                     w_{0} 	& \ldots & w_{-j} 	& \ldots	& w_{-n} \cr
                     \vdots	& \vdots & \vdots & \vdots 	& \vdots \cr
                     w_{n-1}	& \ldots & w_{n-j-1} 	& \ldots 	& w_{-1} \cr
                     1 	& \ldots & z^j 	& \ldots 	& z^n  \cr
              \end{pmatrix} ,
  \label{ops_DetRep:a} \\
  & = (-1)^n\kappa_n\frac{I^{0}_{n}[w(\zeta)(\zeta-z)]}{I^{0}_{n}[w(\zeta)]} ,
  \label{ops_IntRep:a}
\end{align}
\begin{align}
  \phi^*_{n}(z) & = \frac{\kappa_n}{I^{0}_{n}}
         \det \begin{pmatrix}
                     w_{0} 	& \ldots & w_{-n+1}	& z^n  \cr
                     \vdots	& \vdots & \vdots 	& \vdots \cr
                     w_{n-j} 	& \ldots & w_{-j+1} 	& z^j  \cr
                     \vdots	& \vdots & \vdots 	& \vdots \cr
                     w_{n} 	& \ldots & w_{1} 	& 1  \cr
              \end{pmatrix} ,
  \label{ops_DetRep:b} \\
  & = \kappa_n\frac{I^{0}_{n}[w(\zeta)(1-z\zeta^{-1})]}{I^{0}_{n}[w(\zeta)]} .
  \label{ops_IntRep:b} 
\end{align}

The system is alternatively defined by the sequence of ratios $ r_n = \phi_n(0)/\kappa_n $,
known as reflection coefficients because of their role in the scattering theory
formulation of OPS on the unit circle, together with a companion quantity
$ \bar{r}_n = \bar{\phi}_n(0)/\kappa_n $. As in the Szeg\"o theory 
\cite{ops_Sz} $ r_n $ and $ \bar{r}_n $ are related to the above Toeplitz 
determinants by
\begin{equation}
  r_n = (-1)^n\frac{I^{1}_{n}[w]}{I^{0}_{n}[w]}, \quad
  \bar{r}_n = (-1)^n\frac{I^{-1}_{n}[w]}{I^{0}_{n}[w]}.
\end{equation}
The Toeplitz determinants of central interest can then be recovered through the 
following result.
\begin{proposition}[\cite{Ba_1961}]
With the convention $ I^{0}_{0}=1 $ the sequence of $ \{I^{0}_{n}\}^{\infty}_{n=0} $ 
satisfy 
\begin{equation}
   \frac{I^{0}_{n+1}[w]I^{0}_{n-1}[w]}{(I^{0}_{n}[w])^2}
   = 1 - r_{n}\bar{r}_n, \quad n \geq 1.
\label{ops_I0}
\end{equation}
subject to the condition $ r_{n}\bar{r}_n \neq 1 $ for all $ n \geq 1 $. 
\end{proposition}

Fundamental consequences of the orthogonality condition are the following coupled linear 
recurrence relations.
\begin{proposition}[\cite{Ba_1960,Ba_1961}]
\begin{align}
  \kappa_n  \phi_{n+1}(z)
   & = \kappa_{n+1}z \phi_{n}(z)+\phi_{n+1}(0) \phi^*_n(z) ,
  \label{ops_rr:a} \\
  \kappa_n \phi^*_{n+1}(z)
   & = \kappa_{n+1} \phi^*_{n}(z)+\bar{\phi}_{n+1}(0) z\phi _n(z) .
  \label{ops_rr:b}
\end{align}
\end{proposition}
One finds three-term or second order recurrences for the uncoupled recurrence relations  
\begin{align}
  \kappa_n\phi_n(0)\phi_{n+1}(z) + \kappa_{n-1}\phi_{n+1}(0)z\phi_{n-1}(z)
   & = [\kappa_{n}\phi_{n+1}(0)+\kappa_{n+1}\phi_{n}(0)z]\phi_n(z) ,
  \label{ops_ttr:a} \\
  \kappa_n\bar{\phi}_n(0)\phi^*_{n+1}(z) + \kappa_{n-1}\bar{\phi}_{n+1}(0)z\phi^*_{n-1}(z)
   & = [\kappa_{n}\bar{\phi}_{n+1}(0)z+\kappa_{n+1}\bar{\phi}_{n}(0)]\phi^*_n(z) .
  \label{ops_ttr:b} 
\end{align}

The analogue of the Christoffel-Darboux summation formula is given by the following
result.
\begin{proposition}[\cite{Ba_1961}]
\begin{align}
  \sum^{n}_{j=0} \phi_j(z)\bar{\phi}_j(\bar{\zeta}) 
  & =
  \frac{\phi^*_n(z)\overline{\phi^*_n}(\bar{\zeta})
        -z\bar{\zeta}\phi_n(z)\bar{\phi}_n(\bar{\zeta})}{1-z\bar{\zeta}} ,
  \label{ops_CD:a} \\
  & = 
  \frac{\phi^*_{n+1}(z)\overline{\phi^*_{n+1}}(\bar{\zeta})
        -\phi_{n+1}(z)\bar{\phi}_{n+1}(\bar{\zeta})}{1-z\bar{\zeta}} ,
  \label{ops_CD:b}
\end{align}
for $ z\bar{\zeta} \not= 1 $.
Here 
\begin{equation}
   \overline{\phi^*_n}(\bar{\zeta}) = \bar{\zeta}^n\phi_n(1/\bar{\zeta}) .
\end{equation}
\end{proposition}

Equations (\ref{ops_ttr:a},\ref{ops_ttr:b}) being second order linear difference 
equations admit other linearly independent solutions $ \psi_n(z), \psi^*_n(z) $, 
and we define two such polynomial solutions, the polynomials of the second kind 
or associated polynomials
\begin{equation}
  \psi_n(z) 
  := \int_{\TT}\frac{d\zeta}{2\pi i\zeta}\frac{\zeta+z}{\zeta-z}w(\zeta)
               [\phi_n(\zeta)-\phi_n(z)],
       \quad n \geq 1, \quad \psi_0 := \kappa_0w_0 = 1/\kappa_0,
\label{ops_psi:a}
\end{equation}
and its reciprocal polynomial $ \psi^*_n(z) $. 
The integral formula for $ \psi^*_{n} $ is
\begin{equation}
  \psi^*_n(z) 
  := -\int_{\TT}\frac{d\zeta}{2\pi i\zeta}\frac{\zeta+z}{\zeta-z}w(\zeta)
                [z^n\bar{\phi}_n(\bar{\zeta})-\phi^*_n(z)],
       \quad n \geq 1, \quad \psi^*_0 := 1/\kappa_0.
\label{ops_psi:b}
\end{equation} 
A central object in the theory is the Carath\'eodory function, or generating 
function of the Toeplitz elements
\begin{equation}
   F(z) := \int_{\TT}\frac{d\zeta}{2\pi i\zeta}\frac{\zeta+z}{\zeta-z}w(\zeta) ,
\label{ops_Cfun}
\end{equation}
which has the expansions inside and outside the unit circle
\begin{equation}
   F(z) = \begin{cases}
            w_0 + 2 \sum^{\infty}_{k=1}w_{k} z^k, 
           & \text{if $ |z| < \Delta_{min} < 1 $}, \\
           -w_0 - 2 \sum^{\infty}_{k=1}w_{-k} z^{-k}, 
           & \text{if $ |z| > \Delta_{max} > 1 $}.
          \end{cases}
\end{equation}
Having these definitions one requires two non-polynomial solutions 
$ \epsilon_n(z), \epsilon^*_n(z) $ to the recurrences and these are constructed 
as linear combinations of the polynomial solutions according to 
\begin{align}
   \epsilon_n(z) := \psi_n(z)+F(z)\phi_n(z)
   & = \int_{\TT}\frac{d\zeta}{2\pi i\zeta}\frac{\zeta+z}{\zeta-z}w(\zeta)
                   \phi_n(\zeta) ,
   \label{ops_eps:a} \\
   \epsilon^*_n(z) := \psi^*_n(z)-F(z)\phi^*_n(z)
   & = -z^n\int_{\TT}\frac{d\zeta}{2\pi i\zeta}\frac{\zeta+z}{\zeta-z}w(\zeta)
                   \bar{\phi}_n(\bar{\zeta}) ,
   \nonumber \\
   & = \frac{1}{\kappa_n} 
          -\int_{\TT}\frac{d\zeta}{2\pi i\zeta}\frac{\zeta+z}{\zeta-z}w(\zeta)
                   \phi^*_n(\zeta) .
   \label{ops_eps:b}
\end{align}
These have integral representations analogous to (\ref{ops_IntRep:a},\ref{ops_IntRep:b}) 
\begin{align}
  \frac{\kappa_n}{2}\epsilon_{n}(z)
  & = z^n\frac{I^{1}_{n+1}[w(\zeta)(\zeta-z)^{-1}]}{I^{0}_{n+1}[w(\zeta)]} ,
  \label{ops_IntRep:c} \\
  \frac{\kappa_n}{2}\epsilon^*_{n}(z)
  & = (-z)^{n+1}\frac{I^{0}_{n+1}[w(\zeta)(\zeta-z)^{-1}]}{I^{0}_{n+1}[w(\zeta)]} .
  \label{ops_IntRep:d} 
\end{align}
for $ |z| \neq 1 $, 
which are particular cases of more general moments of the characteristic polynomial 
considered by Ismail and R\"udemann \cite{IR_1992}.

\begin{theorem}[\cite{Ge_1961},\cite{Ge_1962},\cite{Ge_1977},\cite{JNT_1989}]
$ \psi_n(z), \psi^*_{n}(z)$ satisfy the three-term recurrence relations 
(\ref{ops_ttr:a}, \ref{ops_ttr:b}) and along with 
$ \epsilon_{n}(z), \epsilon^*_{n}(z) $ satisfy a variant of 
(\ref{ops_rr:a},\ref{ops_rr:b}) namely
\begin{align}
  \kappa_n  \epsilon_{n+1}(z)
   & = \kappa_{n+1}z \epsilon_{n}(z)-\phi_{n+1}(0) \epsilon^*_n(z) ,
  \label{ops_rre:a} \\
  \kappa_n \epsilon^*_{n+1}(z)
   & = \kappa_{n+1} \epsilon^*_{n}(z)-\bar{\phi}_{n+1}(0) z\epsilon _n(z) .
  \label{ops_rre:b}
\end{align}
\end{theorem}

\begin{theorem}[\cite{Ge_1961}]
The Casoratians of the polynomial solutions 
$ \phi_n, \phi^*_n, \psi_n, \psi^*_{n} $ are
\begin{align}
   \phi_{n+1}(z)\psi_n(z) - \psi_{n+1}(z)\phi_n(z) 
  & = \phi_{n+1}(z)\epsilon_n(z) - \epsilon_{n+1}(z)\phi_n(z)
    = 2\frac{\phi_{n+1}(0)}{\kappa_n}z^n ,
  \label{ops_Cas:a} \\
   \phi^*_{n+1}(z)\psi^*_n(z) - \psi^*_{n+1}(z)\phi^*_n(z) 
  & = \phi^*_{n+1}(z)\epsilon^*_n(z) - \epsilon^*_{n+1}(z)\phi^*_n(z) 
    = 2\frac{\bar{\phi}_{n+1}(0)}{\kappa_n}z^{n+1} , 
  \label{ops_Cas:b} \\
   \phi_{n}(z)\psi^*_n(z) + \psi_{n}(z)\phi^*_n(z) 
  & = \phi_{n}(z)\epsilon^*_n(z) + \epsilon_{n}(z)\phi^*_n(z)
    = 2z^n .  
  \label{ops_Cas:c}
\end{align} 
\end{theorem}

Further identities from the Szeg\"o theory that generalise are those that relate 
the leading coefficients back to the reflection coefficients
\begin{align}
   \kappa_n^2 & = \kappa_{n-1}^2 + \phi_n(0)\bar{\phi}_n(0),
   \label{ops_kappa}\\
   \frac{l_{n}}{\kappa_{n}} & = \frac{l_{n-1}}{\kappa_{n-1}}+r_{n}\bar{r}_{n-1},
   \label{ops_l}\\
   \frac{m_{n}}{\kappa_{n}} & = \frac{m_{n-1}}{\kappa_{n-1}}
          +r_{n}\Big[ \bar{r}_{n-2} + \bar{r}_{n-1}\frac{l_{n-2}}{\kappa_{n-2}} \Big].
   \label{ops_m}
\end{align}
Some useful relations for the leading coefficients of the product of a monomial
and an bi-orthogonal polynomial or its derivative are 
\begin{align}
\begin{split}
  z\phi_n(z) & = \frac{\kappa_{n}}{\kappa_{n+1}} \phi_{n+1}(z) 
       + \left(\frac{l_{n}}{\kappa_{n}}-\frac{l_{n+1}}{\kappa_{n+1}}\right) \phi_{n}(z)
  \\
  & \quad
       + \left\{ \frac{l_{n}}{\kappa_{n-1}}
                 \left( \frac{l_{n+1}}{\kappa_{n+1}}-\frac{l_{n}}{\kappa_{n}}\right)
                + \frac{m_{n}}{\kappa_{n-1}}-\frac{m_{n+1}}{\kappa_{n+1}}
                  \frac{\kappa_{n}}{\kappa_{n-1}} \right\} \phi_{n-1}(z)
       + \pi_{n-2} ,
  \\
  z^2\phi_n(z) & = \frac{\kappa_{n}}{\kappa_{n+2}} \phi_{n+2}(z) 
       + \left(\frac{l_{n}}{\kappa_{n+1}}-\frac{l_{n+2}}{\kappa_{n+2}}
               \frac{\kappa_{n}}{\kappa_{n+1}}\right) \phi_{n+1}(z)
  \\
  & \quad
       + \left\{ \frac{l_{n+1}}{\kappa_{n+1}}
                 \left( \frac{l_{n+2}}{\kappa_{n+2}}-\frac{l_{n}}{\kappa_{n}}\right)
                + \frac{m_{n}}{\kappa_{n}}-\frac{m_{n+2}}{\kappa_{n+2}}
                  \right\} \phi_{n}(z) + \pi_{n-1} , 
\end{split}
\end{align}
\begin{align}
\begin{split}
  \phi'_n(z) 
  & = n\frac{\kappa_{n}}{\kappa_{n-1}} \phi_{n-1}(z) + \pi_{n-2} , 
  \\
  z\phi'_n(z)
  & = n \phi_n(z) - \frac{l_{n}}{\kappa_{n-1}} \phi_{n-1}(z)
                  + \pi_{n-2} ,
  \\
  z^2\phi'_n(z)
  & = n\frac{\kappa_{n}}{\kappa_{n+1}} \phi_{n+1}(z)
         + \left\{(n-1)\frac{l_{n}}{\kappa_{n}}-n\frac{l_{n+1}}{\kappa_{n+1}}
           \right\} \phi_n(z) + \pi_{n-1} ,
\end{split}
\label{ops_prod}
\end{align}
where $ ' $ denotes the derivative with respect to $ z $ and
where $ \pi_{n} $ denotes an arbitrary polynomial of the linear space of polynomials 
with degree at most $ n $. These identities can be verified directly.
 
We will require the leading order terms in expansions of 
$ \phi_n(z), \phi^*_n(z), \epsilon_n(z), \epsilon^*_{n}(z) $ both inside and 
outside the unit circle.
\begin{corollary}
The bi-orthogonal polynomials $ \phi_n(z), \phi^*_n(z) $ have the following 
expansions
\begin{align}
   \phi_n(z) & = 
   \begin{cases}
      \phi_n(0)
       + \dfrac{1}{\kappa_{n-1}}
         (\kappa_n\phi_{n-1}(0)+\phi_n(0)\bar{l}_{n-1})z
       + {\rm O}(z^2) & |z| < 1 ,\\
      \kappa_n z^n + l_n z^{n-1} + {\rm O}(z^{n-2}) & |z| > 1 ,
   \end{cases} \label{ops_phiexp:a} \\
   \phi^*_n(z) & = 
   \begin{cases}
      \kappa_n + \bar{l}_n z + {\rm O}(z^{2}) & |z| < 1 ,\\
      \bar{\phi}_n(0) z^n
       + \dfrac{1}{\kappa_{n-1}}
         (\kappa_n\bar{\phi}_{n-1}(0)+\bar{\phi}_n(0)l_{n-1})z^{n-1}
       + {\rm O}(z^{n-2}) & |z| > 1 ,\\
   \end{cases} \label{ops_phiexp:b}
\end{align}
whilst the associated functions have the expansions
\begin{align}
   \dfrac{\kappa_n}{2}\epsilon_n(z) & = 
   \begin{cases}
      z^n
       - \dfrac{\bar{l}_{n+1}}{\kappa_{n+1}}z^{n+1}
       + {\rm O}(z^{n+2}) & |z| < 1 ,\\
         \dfrac{\phi_{n+1}(0)}{\kappa_{n+1}}z^{-1}
       + \left(\dfrac{\kappa^2_n}{\kappa^2_{n+1}}
               \dfrac{\phi_{n+2}(0)}{\kappa_{n+2}}
               -\dfrac{\phi_{n+1}(0)}{\kappa_{n+1}}
                \dfrac{l_{n+1}}{\kappa_{n+1}}
         \right)z^{-2} & \\
       \hfill + {\rm O}(z^{-3}) & |z| > 1 ,
   \end{cases} \label{ops_epsexp:a}
\end{align}
\begin{align}
   \dfrac{\kappa_n}{2}\epsilon^*_n(z) & = 
   \begin{cases}
         \dfrac{\bar{\phi}_{n+1}(0)}{\kappa_{n+1}}z^{n+1}
       + \left(\dfrac{\kappa^2_n}{\kappa^2_{n+1}}
                \dfrac{\bar{\phi}_{n+2}(0)}{\kappa_{n+2}}
               -\dfrac{\bar{\phi}_{n+1}(0)}{\kappa_{n+1}}
                \dfrac{\bar{l}_{n+1}}{\kappa_{n+1}}
         \right)z^{n+2} & \\
       \hfill + {\rm O}(z^{n+3}) & |z| < 1 ,\\
       1 - \dfrac{l_{n+1}}{\kappa_{n+1}}z^{-1} + \left(
         \dfrac{l_{n+2}l_{n+1}}{\kappa_{n+2}\kappa_{n+1}}
        -\dfrac{m_{n+2}}{\kappa_{n+2}}
         \right)z^{-2} + {\rm O}(z^{-3}) & |z| > 1 ,
   \end{cases} \label{ops_epsexp:b}
\end{align}
\end{corollary}
\begin{proof}
The second line of (\ref{ops_phiexp:a}) and the first line of (\ref{ops_phiexp:b}) 
follow from the definitions. For the remaining lines of these two formulae we use 
\begin{equation}
   \kappa_{n-1}\phi'_{n}(0)=\kappa_n\phi_{n-1}(0)+\phi_n(0)\bar{l}_{n-1} ,
\end{equation}
which results from differentiating (\ref{ops_rr:a}) and setting $ z=0 $.
The first line of (\ref{ops_epsexp:a}) and the second line of (\ref{ops_epsexp:b})
can be derived by employing the uncoupled recurrence relations 
(\ref{ops_ttr:a}) and (\ref{ops_ttr:b}) respectively. The remaining lines of these
two formulae can be found by using the identity
\begin{equation}
   z\phi_n(z) = \frac{\kappa_n}{\kappa_{n+1}}\phi_{n+1}(z)
  - \frac{\phi_{n+1}(0)}{\kappa_n\kappa_{n+1}}\sum^{n}_{j=0}\bar{\phi}_j(0)\phi_j(z) ,
\end{equation}
which in turn follows from combining (\ref{ops_rr:a}) and (\ref{ops_rr:b}).
\end{proof}

The $z$-derivatives or spectral derivatives of the bi-orthogonal polynomials in 
general are related to two consecutive polynomials \cite{IW_2001} and we extend
this to all $\phi_n(z)$, $\phi^*_n(z)$, $\epsilon_n(z)$, $\epsilon^*_n(z)$ with 
the following parameterisation.
\begin{proposition}\label{ops_spectralD}
The derivatives of the bi-orthogonal polynomials and associated functions are
expressible as linear combinations in a related way ($\; ' := d/dz $), 
\begin{align}
   W(z)\phi'_n(z) & = 
   \Theta_n(z) \phi_{n+1}(z) - (\Omega_n(z)+V(z)) \phi_n(z) ,
   \label{ops_zD:a} \\
   W(z)\phi^*_n{\!'}(z) & = 
   -\Theta^*_n(z) \phi^*_{n+1}(z) + (\Omega^*_n(z)-V(z)) \phi^*_n(z) ,
   \label{ops_zD:c} \\
   W(z)\epsilon'_n(z) & = 
   \Theta_n(z) \epsilon_{n+1}(z) - (\Omega_n(z)-V(z)) \epsilon_n(z) ,
   \label{ops_zD:b} \\
   W(z)\epsilon^*_n{\!'}(z) & = 
   -\Theta^*_n(z) \epsilon^*_{n+1}(z) + (\Omega^*_n(z)+V(z)) \epsilon^*_n(z) ,
   \label{ops_zD:d} 
\end{align}
with coefficient functions $ W(z), V(z) $ independent of $ n $.
\end{proposition}
\begin{proof}
The first, (\ref{ops_zD:a}), was found in \cite{IW_2001} where the coefficients
were taken to be (their notation $ A_n, B_n $ should not be confused with our 
use of it subsequently)
\begin{align}
   A_n & = -\frac{\kappa_{n-1}\phi_{n+1}(0)}{\kappa_{n}\phi_{n}(0)}
            \frac{z\Theta_n(z)}{W(z)} ,
  \\
   B_n & = \frac{1}{W(z)}\left( \Omega_n(z)+V(z)
            - \left[\frac{\phi_{n+1}(0)}{\phi_{n}(0)}+\frac{\kappa_{n+1}}{\kappa_{n}}z 
              \right]\Theta_n(z) \right) .
\end{align}
The other differential relations can be found in an analogous manner.
\end{proof}

The coefficient functions 
$ \Theta_n(z), \Omega_n(z), \Theta^*_n(z), \Omega^*_n(z) $
satisfy coupled linear recurrence relations themselves, one of
which was reported in \cite{IW_2001}. The full set are given in the following 
proposition.
\begin{proposition}\label{ops_Linear1}
The coefficient functions satisfy the coupled linear recurrence relations 
\begin{equation}
  \Omega_n(z) + \Omega_{n-1}(z) 
  - \left( \frac{\phi_{n+1}(0)}{\phi_{n}(0)}+\frac{\kappa_{n+1}}{\kappa_{n}}z
    \right)\Theta_n(z) + (n-1)\frac{W(z)}{z} = 0 ,
\label{ops_rrCf:a}
\end{equation}
\begin{multline}
   \left( \frac{\phi_{n+1}(0)}{\phi_{n}(0)}+\frac{\kappa_{n+1}}{\kappa_{n}}z
   \right) (\Omega_{n-1}(z) - \Omega_{n}(z)) 
   \\
   + \frac{\kappa_{n}\phi_{n+2}(0)}{\kappa_{n+1}\phi_{n+1}(0)}z\Theta_{n+1}(z)
   - \frac{\kappa_{n-1}\phi_{n+1}(0)}{\kappa_{n}\phi_{n}(0)}z\Theta_{n-1}(z)
   - \frac{\phi_{n+1}(0)}{\phi_{n}(0)}\frac{W(z)}{z} = 0 ,
\label{ops_rrCf:b}
\end{multline}
\begin{equation}
  \Omega^*_n(z) + \Omega^*_{n-1}(z) 
  - \left( \frac{\kappa_{n+1}}{\kappa_{n}}+\frac{\bar{\phi}_{n+1}(0)}{\bar{\phi}_{n}(0)}z
    \right)\Theta^*_n(z) - n\frac{W(z)}{z} = 0 ,
\label{ops_rrCf:c}
\end{equation}
\begin{multline}
   \left( \frac{\kappa_{n+1}}{\kappa_{n}}+\frac{\bar{\phi}_{n+1}(0)}{\bar{\phi}_{n}(0)}z
   \right) (\Omega^*_{n-1}(z) - \Omega^*_{n}(z)) 
   \\
   + \frac{\kappa_{n}\bar{\phi}_{n+2}(0)}{\kappa_{n+1}\bar{\phi}_{n+1}(0)}
      z\Theta^*_{n+1}(z)
   - \frac{\kappa_{n-1}\bar{\phi}_{n+1}(0)}{\kappa_{n}\bar{\phi}_{n}(0)}
      z\Theta^*_{n-1}(z)
   + \frac{\kappa_{n+1}}{\kappa_{n}}\frac{W(z)}{z} = 0 ,
\label{ops_rrCf:d}
\end{multline}
\begin{equation}
  \Omega_{n+1}(z) + \Omega^*_{n}(z) 
  - \left( \frac{\phi_{n+2}(0)}{\phi_{n+1}(0)}+\frac{\kappa_{n+2}}{\kappa_{n+1}}z
    \right)\Theta_{n+1}(z)
  + \frac{\kappa_{n+1}}{\kappa_{n}}(z\Theta_n(z)-\Theta^*_n(z)) = 0 ,
\label{ops_rrCf:e}
\end{equation}
\begin{multline}
   \Omega_{n}(z) - \Omega_{n+1}(z) 
   + \frac{\kappa_{n+2}}{\kappa_{n+1}}
    \left( z+\frac{\bar{\phi}_{n+1}(0)}{\kappa_{n+1}}\frac{\phi_{n+2}(0)}{\kappa_{n+2}}
    \right)\Theta_{n+1}(z) \\
   + \frac{\phi_{n+1}(0)\bar{\phi}_{n+1}(0)}{\kappa_{n+1}\kappa_{n}}\Theta^*_n(z)
   - \frac{\kappa_{n+1}}{\kappa_{n}}z\Theta_{n}(z)
   - \frac{W(z)}{z} = 0 ,
\label{ops_rrCf:f}
\end{multline}
\begin{multline}
  \Omega^*_{n+1}(z) + \Omega_{n}(z) 
  - \left( \frac{\kappa_{n+2}}{\kappa_{n+1}}
           +\frac{\bar{\phi}_{n+2}(0)}{\bar{\phi}_{n+1}(0)}z
    \right)\Theta^*_{n+1}(z) \\
  - \frac{\kappa_{n+1}}{\kappa_{n}}(z\Theta_n(z)-\Theta^*_n(z)) - \frac{W(z)}{z}= 0 ,
\label{ops_rrCf:g}
\end{multline}
\begin{multline}
   \Omega^*_{n}(z) - \Omega^*_{n+1}(z)
   + \frac{\kappa_{n+2}}{\kappa_{n+1}}
    \left( 1+\frac{\phi_{n+1}(0)}{\kappa_{n+1}}\frac{\bar{\phi}_{n+2}(0)}{\kappa_{n+2}}z
    \right)\Theta^*_{n+1}(z) \\
   + \frac{\phi_{n+1}(0)\bar{\phi}_{n+1}(0)}{\kappa_{n+1}\kappa_{n}}z\Theta_n(z)
   - \frac{\kappa_{n+1}}{\kappa_{n}}\Theta^*_{n}(z) = 0 .
\label{ops_rrCf:h}
\end{multline}
\end{proposition} 
\begin{proof}
The first (\ref{ops_rrCf:a}) was found in \cite{IW_2001} by a direct evaluation
of the left-hand side using integral definitions of the coefficient functions,
however all of the relations follow from the compatibility of the differential
relations and the recurrence relations. 
Thus (\ref{ops_rrCf:a},\ref{ops_rrCf:b}) follow from the compatibility of 
(\ref{ops_zD:a}) and (\ref{ops_ttr:a}), 
(\ref{ops_rrCf:c}),(\ref{ops_rrCf:d}) follow from 
(\ref{ops_zD:c}) and (\ref{ops_ttr:b}),
(\ref{ops_rrCf:e},\ref{ops_rrCf:f}) follow from the combination of
(\ref{ops_zD:a},\ref{ops_zD:c}) and (\ref{ops_rr:a}), and 
(\ref{ops_rrCf:g},\ref{ops_rrCf:h}) follow from the combination of
(\ref{ops_zD:a},\ref{ops_zD:c}) and (\ref{ops_rr:b}).
\end{proof}

\begin{remark}
The relations given above are obviously not all independent, as for example we
note that (\ref{ops_rrCf:a}) can derived from (\ref{ops_rrCf:e}) with the use
of (\ref{ops_rrCf:j}) below.
\end{remark}

\begin{corollary}\label{ops_Linear2}
Some additional identities satisfied by the coefficient functions are the following
\begin{gather}
  \frac{\phi_{n+1}(0)}{\phi_{n}(0)}\Theta_{n}(z)
 - \frac{\kappa_{n}}{\kappa_{n-1}}z\Theta_{n-1}(z) 
 = \frac{\bar{\phi}_{n+1}(0)}{\bar{\phi}_{n}(0)}z\Theta^*_n(z)
  - \frac{\kappa_{n}}{\kappa_{n-1}} \Theta^*_{n-1}(z) ,
 \label{ops_rrCf:i} \\  
  \Omega^*_{n}(z)-\Omega_{n}(z)
  = - \frac{\kappa_{n+1}}{\kappa_{n}}(z\Theta_{n}(z)-\Theta^*_n(z))+n\frac{W(z)}{z} ,
 \label{ops_rrCf:j} \\
  \Omega^*_{n}(z)+\Omega_{n}(z)  
  = \frac{\kappa^2_{n}}{\kappa^2_{n+1}}
    \left[\frac{\phi_{n+2}(0)}{\phi_{n+1}(0)}\Theta_{n+1}(z)
          + \frac{\kappa_{n+1}}{\kappa_{n}} \Theta^*_{n}(z)\right] 
    + \frac{W(z)}{z} .
 \label{ops_rrCf:k}
\end{gather}
\end{corollary}
\begin{proof}
Consider (\ref{ops_rrCf:i}) first. In (\ref{ops_rrCf:a}) let us map $ n \mapsto n+1 $
and add this result to (\ref{ops_rrCf:c}). However the combination
$ \Omega_{n+1}(z)+\Omega^*_{n}(z) $, which appears in the resulting sum in this form and
also under $ n \mapsto n-1 $, occurs in (\ref{ops_rrCf:e}) allowing it to be 
eliminated. To derive (\ref{ops_rrCf:j}) we again map $ n \mapsto n+1 $ in (\ref{ops_rrCf:a})
and subtract this from (\ref{ops_rrCf:e}) eliminating $ \Omega_{n+1}(z) $.
If we add (\ref{ops_rrCf:e}) and (\ref{ops_rrCf:f}) and simplify the result then
we arrive at (\ref{ops_rrCf:k}).
\end{proof}

For a general system of bi-orthogonal polynomials on the unit circle the coupled 
recurrence relations and spectral differential relations 
can be reformulated in terms of first order $2\times 2$ matrix equations (or
alternatively as second order scalar equations). Here we define our
matrix variables and derive such matrix relations, and this serves as an 
introduction to a characterisation of the general bi-orthogonal polynomial system on
the unit circle as the solution to a $2\times 2$ matrix Riemann-Hilbert problem.

Firstly we note that the recurrence relations for the associated functions 
$ \epsilon_n(z) $, $ \epsilon^*_n(z) $ given in (\ref{ops_rre:a},\ref{ops_rre:b}) differ
from those of the polynomial systems (\ref{ops_rr:a},\ref{ops_rr:b}) by a reversal of
the signs of $ \phi_n(0) $, $ \bar{\phi}_n(0) $. We can compensate for this by constructing
the $ 2\times 2 $ matrix
\begin{equation}
   Y_n(z) := 
   \begin{pmatrix}
          \phi_n(z)   & \dfrac{\epsilon_n(z)}{w(z)} \cr
          \phi^*_n(z) & -\dfrac{\epsilon^{\vphantom{I}*}_n(z)}{w(z)} \cr
   \end{pmatrix} ,
\label{ops_Ydefn}
\end{equation}
and note from (\ref{ops_Cas:b}) that $ \det Y_n = -2z^n/w(z) $.

\begin{corollary}
The recurrence relations for a general system of bi-orthogonal polynomials 
(\ref{ops_rr:a},\ref{ops_rr:b}) and their associated functions 
(\ref{ops_rre:a},\ref{ops_rre:b}) are equivalent to the matrix recurrence
\begin{equation}
   Y_{n+1} := K_n Y_{n}
   = \frac{1}{\kappa_n}
       \begin{pmatrix}
              \kappa_{n+1} z   & \phi_{n+1}(0) \cr
              \bar{\phi}_{n+1}(0) z & \kappa_{n+1} \cr
       \end{pmatrix} Y_{n} . 
\label{ops_Yrecur}
\end{equation} 
According to (\ref{ops_kappa}) the matrix $ K_n $ has the property $ \det K_n = z $.
\end{corollary}

\begin{corollary}
The system of spectral derivatives for a general system of bi-orthogonal polynomials
and associated functions (\ref{ops_zD:a}-\ref{ops_zD:d}) are equivalent to 
the matrix differential equation
\begin{multline}
   Y'_{n} := A_n Y_{n} \\
   = \frac{1}{W(z)}
       \begin{pmatrix}
              -\left[ \Omega_n(z)+V(z)
                     -\dfrac{\kappa_{n+1}}{\kappa_n}z\Theta_n(z)
               \right]
            & \dfrac{\phi_{n+1}(0)}{\kappa_n}\Theta_n(z)
            \cr
              -\dfrac{\bar{\phi}_{n+1}(0)}{\kappa_n}z\Theta^*_n(z)
            &  \Omega^*_n(z)-V(z)
                     -\dfrac{\kappa_{n+1}}{\kappa_n}\Theta^*_n(z)
            \cr
       \end{pmatrix} Y_{n} . 
\label{ops_YzDer}
\end{multline} 
\end{corollary}
\begin{proof}
This follows from (\ref{ops_zD:a}-\ref{ops_zD:d}) and employing 
(\ref{ops_rr:a},\ref{ops_rr:b},\ref{ops_rre:a},\ref{ops_rre:b}).
We note that $ {\rm Tr}A_n = n/z-w'/w $ when (\ref{ops_rrCf:j}) is employed.
\end{proof}

\begin{remark}
Compatibility of the relations (\ref{ops_Yrecur}) and (\ref{ops_YzDer}) leads to 
\begin{equation}
   K'_n = A_{n+1}K_n-K_nA_n ,
\end{equation}
and upon examining the $ (1,1) $-component of this we recover the linear recurrence 
(\ref{ops_rrCf:f}), the $ (1,2) $-component yields (\ref{ops_rrCf:e}), whilst the
$ (2,1) $-component gives (\ref{ops_rrCf:g}) and the $ (2,2) $-component implies 
(\ref{ops_rrCf:h}).
\end{remark}

\begin{remark}
There are, in a second-order difference equation such as (\ref{ops_rre:a}) or 
(\ref{ops_rre:b}), other forms of the matrix variables and equations and these
alternative forms will appear in our subsequent work. Defining
\begin{equation}
   X_n(z;t) := 
   \begin{pmatrix}
          \phi_{n+1}(z)   & \dfrac{\epsilon_{n+1}(z)}{w(z)} \cr
          \phi_n(z) & \dfrac{\epsilon_n(z)}{w(z)} \cr
   \end{pmatrix} ,
   \quad
   X^*_n(z;t) := 
   \begin{pmatrix}
          \phi^*_{n+1}(z)   & \dfrac{\epsilon^*_{n+1}(z)}{w(z)} \cr
          \phi^*_n(z) & \dfrac{\epsilon^*_n(z)}{w(z)} \cr
   \end{pmatrix} ,
\label{ops_Xdefn}
\end{equation}
we find the spectral derivatives to be
\begin{equation}
   W(z)X'_{n} = \begin{pmatrix}
              \Omega_n(z)-V(z)
                     +n\dfrac{W(z)}{z}
            & -\dfrac{\kappa_n\phi_{n+2}(0)}{\kappa_{n+1}\phi_{n+1}(0)} z\Theta_{n+1}(z)
            \cr
              \Theta_n(z)
            & -\Omega_n(z)-V(z)
            \cr
                \end{pmatrix} X_{n} , 
\label{ops_XzDer:a}
\end{equation} 
\begin{equation}
   W(z)X^*_{n}{\!'} = \begin{pmatrix}
              -\Omega^*_n(z)-V(z)
                     +(n+1)\dfrac{W(z)}{z}
            & \dfrac{\kappa_n\bar{\phi}_{n+2}(0)}{\kappa_{n+1}\bar{\phi}_{n+1}(0)}
              z\Theta^*_{n+1}(z)
            \cr
              -\Theta^*_n(z)
            & \Omega^*_n(z)-V(z)
            \cr
                      \end{pmatrix} X^*_{n} . 
\label{ops_XzDer:b}
\end{equation} 
Another system is based upon the definition
\begin{equation}
   Z_n(z;t) := 
   \begin{pmatrix}
          \phi_{n+1}(z)   & \dfrac{\epsilon_{n+1}(z)}{w(z)} \cr
          \phi^*_n(z) & -\dfrac{\epsilon^*_n(z)}{w(z)} \cr
   \end{pmatrix} ,
   \quad
   Z^*_n(z;t) := 
   \begin{pmatrix}
          \phi^*_{n+1}(z)   & -\dfrac{\epsilon^*_{n+1}(z)}{w(z)} \cr
          \phi_n(z) & \dfrac{\epsilon_n(z)}{w(z)} \cr
   \end{pmatrix} ,
\label{ops_Zdefn}
\end{equation}
and in this case the spectral derivatives are
\begin{multline}
   W(z)Z'_{n} = \\
       \begin{pmatrix}
              -\Omega^*_n(z)-V(z)
              +\dfrac{\kappa_n}{\kappa_{n+1}}\Theta^*_n(z)
                     +(n+1)\dfrac{W(z)}{z}
            & \dfrac{\kappa_n\phi_{n+2}(0)}{\kappa^2_{n+1}}\Theta_{n+1}(z)
            \cr
              -\dfrac{\bar{\phi}_{n+1}(0)}{\kappa_{n+1}}\Theta^*_n(z)
            & \Omega^*_n(z)-V(z)
              -\dfrac{\kappa_n}{\kappa_{n+1}}\Theta^*_n(z)
            \cr
       \end{pmatrix} \\ \times Z_{n} , 
\label{ops_ZzDer:a}
\end{multline} 
\begin{multline}
   W(z)Z^*_{n}{\!'} \\
     = \begin{pmatrix}
              \Omega_n(z)-V(z)
              -\dfrac{\kappa_n}{\kappa_{n+1}}z\Theta_n(z)
            & -\dfrac{\kappa_n\bar{\phi}_{n+2}(0)}{\kappa^2_{n+1}}z^2\Theta^*_{n+1}(z)
            \cr
              \dfrac{\phi_{n+1}(0)}{\kappa_{n+1}}\Theta_n(z)
            & -\Omega_n(z)-V(z)
              +\dfrac{\kappa_n}{\kappa_{n+1}}z\Theta_n(z)
            \cr
       \end{pmatrix} Z^*_{n} . 
\label{ops_ZzDer:b}
\end{multline} 
\end{remark}

We end this section with a characterisation of a general system of bi-orthogonal 
polynomials on the unit circle (and their associated functions) as a solution to 
a particular Riemann-Hilbert problem.

\begin{proposition}\label{ops_RHP}
Consider the following Riemann-Hilbert problem for a $ 2 \times 2 $ matrix function
$ Y : {\CC} \to SL(2,\CC) $ defined in the following statements
\begin{enumerate}
  \item
   $ Y(z) $ is analytic in $ \{z: |z| > 1\}\cup\{z: |z|< 1\} $,
  \item
   on $ z \in \Sigma $ where $ \Sigma $ is the oriented unit circle in a 
   counter-clockwise sense and $ +(-) $ denote the left(right)-hand side or 
   interior(exterior) 
   \begin{equation}
      Y_{+}(z) = Y_{-}(z)
       \begin{pmatrix}
              1 & w(z)/z \cr
              0 & 1 \cr
       \end{pmatrix} ,
   \label{ops_RHjump}
   \end{equation}
  \item
   as $ z \to \infty $
   \begin{equation}
     Y(z) = \left( \mathbb{I} + {\rm O}(z^{-1}) \right)
       \begin{pmatrix}
              z^n             & {\rm O}(z^{-2}) \cr
              {\rm O}(z^{n}) & -z^{-1}\cr
       \end{pmatrix} ,
   \label{ops_RHzLarge}
   \end{equation}
  \item
   as $ z \to 0 $
   \begin{equation}
     Y(z) = \left( \mathbb{I} + {\rm O}(z) \right)
       \begin{pmatrix}
              {\rm O}(1) & {\rm O}(z^{n-1}) \cr
              {\rm O}(1) & {\rm O}(z^{n})   \cr
       \end{pmatrix} .
   \label{ops_RHzSmall}
   \end{equation}
\end{enumerate}
It is assumed that the weight function $ w(z) $ satisfies the restrictions given
at the beginning of this section. Then the unique solution to this Riemann-Hilbert 
problem is 
\begin{equation}
   Y(z) = 
       \begin{pmatrix}
               \dfrac{\phi_n(z)}{\kappa_n}
            &  \dfrac{\epsilon_n(z)}{2\kappa_n z} \\
               \kappa_n\phi^*_n(z) 
            & -\dfrac{\kappa_n\epsilon^{\vphantom{I}*}_n(z)}{2z}
       \end{pmatrix} , \quad n \geq 1 .
\end{equation}
\end{proposition}
\begin{proof}
We firstly note from the jump condition (\ref{ops_RHjump}) that $ Y_{11}, Y_{21} $ 
are entire $ z \in \CC $. From the $ (1,1) $-entry of the asymptotic condition 
(\ref{ops_RHzLarge}) it is clear that $ Y_{11} = \pi_n(z) $ is a polynomial of degree
at most $ n $. Similarly $ Y_{21} = \sigma_n(z) $ from an observation of the
$ (2,1) $-component. From the $ (1,2) $- and $ (2,2) $-components of the jump condition
we deduce 
\begin{equation}
   Y_{+12}-Y_{-12} = \frac{w(z)}{z}Y_{11} ,
  \qquad
   Y_{+22}-Y_{-22} = \frac{w(z)}{z}Y_{21} ,
\end{equation}
and therefore
\begin{equation}
   Y_{12} 
  = \int_{\TT}\frac{d\zeta}{2\pi i\zeta}\frac{w(\zeta)\pi_n(\zeta)}{\zeta-z} ,
  \qquad
   Y_{22} 
  = \int_{\TT}\frac{d\zeta}{2\pi i\zeta}\frac{w(\zeta)\sigma_n(\zeta)}{\zeta-z} .
\end{equation}
Consider the large $ z $ expansion of $ Y_{12} $ implied by the first of these 
formulae
\begin{equation}
   Y_{12} = -z^{-1}\int_{\TT}\frac{d\zeta}{2\pi i\zeta}w(\zeta)\pi_n(\zeta)
   + {\rm O}(z^{-2}) .
\end{equation}
According to (\ref{ops_RHzLarge}) the integral vanishes and so $ \pi_n(\zeta) $ 
is orthogonal to the monomial $ \zeta^0 $. Now take the small $ z $ expansion
\begin{equation}
   Y_{12} 
   = \sum^{n-2}_{l=0}z^{l}\int_{\TT}\frac{d\zeta}{2\pi i\zeta}
     w(\zeta)\pi_n(\zeta)\overline{\zeta^{l+1}}
   + {\rm O}(z^{n-1}) .
\end{equation}
From the $ (1,2) $-component of the condition (\ref{ops_RHzSmall}) we observe 
that all terms in the sum vanish and we conclude the $ \pi_n(\zeta) $ is 
orthogonal to the monomials $ \zeta, \ldots, \zeta^{n-1} $ and the first 
term which survives has the monomial $ \zeta^{n} $. Thus 
$ \pi_n(\zeta) \propto \phi_n(\zeta) $, and from the explicit coefficient in the 
$ (1,1) $-entry of (\ref{ops_RHzLarge}) $ \pi_n(\zeta) $ is the monic bi-orthogonal polynomial
$ \phi_n(\zeta)/\kappa_n $. We turn our attention to $ Y_{22} $ and examine
the small $ z $ expansion
\begin{equation}
   Y_{22} 
   = \sum^{n-1}_{l=0}z^{l}\int_{\TT}\frac{d\zeta}{2\pi i\zeta}
     w(\zeta)\sigma_n(\zeta)\overline{\zeta^{l+1}}
   + {\rm O}(z^{n}) .
\end{equation}
The $ (2,2) $-component of (\ref{ops_RHzSmall}) tells us that all terms in the
sum vanish and consequently $ \sigma_n(\zeta) $ is orthogonal to all monomials
$ \zeta, \ldots, \zeta^{n} $. Therefore $ \sigma_n(\zeta) \propto \phi^*_n(\zeta) $
and we can determine the proportionality constant from the $ (2,2) $-component of 
the asymptotic formula (\ref{ops_RHzLarge}), by comparing it with
\begin{equation}
   Y_{22} = -z^{-1}\int_{\TT}\frac{d\zeta}{2\pi i\zeta}w(\zeta)\sigma_n(\zeta)
   + {\rm O}(z^{-2}),
\end{equation}
to conclude $ \sigma_n(\zeta) = \kappa_n\phi^*_n(\zeta) $. Finally we note that
\begin{equation}
  \int_{\TT}\frac{d\zeta}{2\pi i\zeta}\frac{w(\zeta)\phi_n(\zeta)}{\zeta-z}
  = \frac{1}{2z}\epsilon_n(z),
  \qquad
  \int_{\TT}\frac{d\zeta}{2\pi i\zeta}\frac{w(\zeta)\phi^*_n(\zeta)}{\zeta-z}
  = -\frac{1}{2z}\epsilon^*_n(z),
\end{equation}
when $ n > 0 $. We also point out $ \det Y(z) = -z^{n-1} $.
\end{proof}

\begin{remark}
Our original matrix solution $ Y_n $ specified by (\ref{ops_Ydefn}) is related 
to the solution of the above Riemann-Hilbert problem by
\begin{equation}
    Y_n(z) = 
       \begin{pmatrix}
              \kappa_n & 0 \cr
              0 & \dfrac{1}{\kappa_n} \cr
       \end{pmatrix}
        Y(z) 
       \begin{pmatrix}
              1 & 0 \cr
              0 & \dfrac{2z}{w(z)} \cr
       \end{pmatrix} , \quad
    Y(z) = 
       \begin{pmatrix}
              \dfrac{1}{\kappa_n} & 0 \cr
              0 & \kappa_n \cr
       \end{pmatrix}
        Y_n(z) 
       \begin{pmatrix}
              1 & 0 \cr
              0 & \dfrac{w(z)}{2z} \cr
       \end{pmatrix} .
\end{equation}
Our formulation of the Riemann-Hilbert problem differs from those given in 
studies concerning orthogonal polynomial systems on the unit circle with
more specialised weights \cite{Ba_2001},\cite{BDMcLMZ_2001},\cite{BDJ_2000a},
\cite{BDJ_1999} and \cite{BKMcLM_2004}. We have chosen this 
formulation as it is closest to that occurring for orthogonal polynomial
systems of the line \cite{De_1999}, the jump matrix is independent of the index $ n $
which only appears in the asymptotic condition and it is simply related to
our matrix formulation (\ref{ops_Ydefn}).
\end{remark}

\section{Regular Semi-classical Weights and Isomonodromic Deformations}
\label{SC+IDsection}
\setcounter{equation}{0}

All of the above results apply for a general class of weights on the unit circle
but now we want to consider an additional restriction, namely the special 
structure of regular or generic semi-classical weights.
\begin{definition}\label{rSC_defn}[\cite{Ma_1995a}]
The log-derivative of a regular or generic semi-classical weight function $ w(z) $
is rational in $ z $ with
\begin{equation}
   W(z)w'(z) = 2V(z)w(z) ,
\label{ops_scwgt}
\end{equation}
where $ V(z), W(z) $ are polynomials with the following properties
\begin{enumerate}
 \item
  $ {\rm deg}(W) \geq 2 $,
 \item
  $ {\rm deg}(V) < {\rm deg}(W) $,
 \item
  the $ m $ zeros of $ W(z) $, $ \{z_1, z_2, \ldots ,z_m\} $ are distinct,
 \item
  the residues $ \rho_k = 2V(z_k)/W'(z_k) \notin \ZZ_{\geq 0} $. 
\end{enumerate}
\end{definition}
The terminology regular refers to the connection of this definition with systems
of linear second order differential equations in the complex plane which possess 
only isolated regular singularities, and we will see the appearance of these
later. An explicit example of such a weight is that of the form
\begin{equation}
  w(z) = \prod^m_{j=1}(z-z_j)^{\rho_j} .
\label{ops_scwgt1}
\end{equation}
Consistent with the requirements of Definition 3.1, we set 
\begin{equation}
   W(z) = \prod^m_{j=1}(z-z_j), \qquad
   \frac{2V(z)}{W(z)} = \sum^m_{j=1}\frac{\rho_j}{z-z_j} .
\label{ops_scwgt2}
\end{equation}
To reduce the computational labour and render the ensuing formulae simpler in
appearance
we are going to assume henceforth that one of the finite singular points is 
fixed at the origin, i.e. that $ W(0) = 0 $.

It follows from these definitions that the Carath\'eodory function satisfies
an inhomogeneous form of (\ref{ops_scwgt}).
\begin{lemma}[\cite{La_1972b},\cite{Ma_1995a}]
Let the weight $ w(z) $ be such that $ w(e^{\pi i})=w(e^{-\pi i}) $. 
The Carath\'eodory function (\ref{ops_Cfun}) satisfies the first order linear
ordinary differential equation 
\begin{equation}
   W(z)F'(z) = 2V(z)F(z)+U(z) ,
\label{ops_FD}
\end{equation}
where $ U(z) $ is a polynomial in $ z $. 
\end{lemma}
\begin{proof}
Following \cite{Ma_1995a} we can write
\begin{equation}
  W(z)F(z) = \int\frac{d\zeta}{2\pi i\zeta}\frac{\zeta+z}{\zeta-z}W(\zeta)w(\zeta)
           + \int\frac{d\zeta}{2\pi i\zeta}\frac{\zeta+z}{\zeta-z}[W(z)-W(\zeta)]w(\zeta) .
\end{equation}
The last term is a polynomial in $ z $ of bounded degree which will be denoted
by $ \pi_3(z) $ as we are not going to evaluate it explicitly. Differentiating
this yields
\begin{equation}
  \frac{d}{dz}W(z)F(z) 
  = \int\frac{d\zeta}{2\pi i\zeta}\frac{2\zeta}{(\zeta-z)^2}W(\zeta)w(\zeta)
           + \pi_2(z) ,
\end{equation}
and we rewrite the first term of the right-hand side as 
\begin{align}
   \int\frac{d\zeta}{2\pi i\zeta}\frac{2\zeta}{(\zeta-z)^2}W(\zeta)w(\zeta)
 & = \frac{1}{2\pi iz}\int d\zeta\frac{2z}{(\zeta-z)^2}W(\zeta)w(\zeta) ,
 \nonumber\\
 & = -\frac{1}{2\pi iz}\int d\zeta\frac{d}{d\zeta}\frac{\zeta+z}{\zeta-z}W(\zeta)w(\zeta) ,
 \nonumber\\
 & = \frac{1}{2\pi iz}\int d\zeta\frac{\zeta+z}{\zeta-z}\frac{d}{d\zeta}(W(\zeta)w(\zeta)) .
\end{align}
Consequently we have
\begin{align}
  \frac{d}{dz}W(z)F(z) 
  & = \frac{1}{z}\int\frac{d\zeta}{2\pi i\zeta}\frac{\zeta+z}{\zeta-z}
                                \zeta\frac{d}{d\zeta}(W(\zeta)w(\zeta)) + \pi_2(z) ,
  \nonumber\\
  & = \frac{1}{z}\int\frac{d\zeta}{2\pi i\zeta}\frac{\zeta+z}{\zeta-z}
                                \zeta[W'(\zeta)+2V(\zeta)]w(\zeta) + \pi_2(z) ,
  \nonumber\\
  & = \frac{1}{z}\int\frac{d\zeta}{2\pi i\zeta}\frac{\zeta+z}{\zeta-z}
                                z[W'(z)+2V(z)]w(\zeta)
  \nonumber\\
  & \phantom{=}\quad
      +\frac{1}{z}\int\frac{d\zeta}{2\pi i\zeta}\frac{\zeta+z}{\zeta-z}
                  \left\{ \zeta[W'(\zeta)+2V(\zeta)]-z[W'(z)+2V(z)] \right\}w(\zeta)
  \nonumber\\
  & \phantom{=}\qquad
      + \pi_2(z) ,
  \nonumber\\
  & = [W'(z)+2V(z)]\int\frac{d\zeta}{2\pi i\zeta}\frac{\zeta+z}{\zeta-z}w(\zeta)
      + \pi_1(z) ,
\end{align}
and (\ref{ops_FD}) follows.
\end{proof}

This lemma leads to the following important set of formulae.
\begin{proposition}\label{ops_SCpoly}
Assume that $ W(0) = 0 $. The coefficient functions 
$ \Theta_n(z) $, $ \Theta^*_n(z) $, $ \Omega_n(z) $, $ \Omega^*_n(z) $ are polynomials in
$ z $ of degree $ m-2, m-2, m-1, m-1 $ respectively. Specifically these have 
leading and trailing expansions of the form
\begin{multline}
  \Theta_n(z) =
 (n+1+\sum^m_{j=1}\rho_j)\frac{\kappa_n}{\kappa_{n+1}}z^{m-2}
 \\
 + \bigg\{ -[(n+1+\sum^m_{j=1}\rho_j)\sum^m_{j=1}z_j - \sum^m_{j=1}\rho_j z_j]
             \frac{\kappa_n}{\kappa_{n+1}}
            +(n+2+\sum^m_{j=1}\rho_j)\frac{\kappa^3_n}{\kappa^2_{n+1}\kappa_{n+2}}
             \frac{\phi_{n+2}(0)}{\phi_{n+1}(0)}
 \\
 -(n+\sum^m_{j=1}\rho_j)\frac{\phi_{n+1}(0)\bar{\phi}_{n}(0)}{\kappa_{n+1}\kappa_{n}}
           -2\frac{\kappa_nl_{n+1}}{\kappa^2_{n+1}}   
   \bigg\}z^{m-3}
 \\ + {\rm O}(z^{m-4}) ,
\label{ops_Thexp:a}
\end{multline}
\begin{multline}
 \Theta_n(z) =
 [2V(0)-nW'(0)]\frac{\phi_{n}(0)}{\phi_{n+1}(0)} 
 \\
 + \bigg\{ [2V'(0)-\shalf nW''(0)]\frac{\phi_{n}(0)}{\phi_{n+1}(0)}
          + [2V(0)-(n-1)W'(0)]\frac{\kappa_n\phi_{n-1}(0)}{\kappa_{n-1}\phi_{n+1}(0)}
 \\
 + \left( [(n+1)W'(0)-2V(0)]\frac{\bar{l}_{n+1}}{\kappa_{n+1}}
                   -[(n-1)W'(0)-2V(0)]\frac{\bar{l}_{n-1}}{\kappa_{n+1}} \right)
           \frac{\phi_{n}(0)}{\phi_{n+1}(0)} 
   \bigg\}z
 \\
 + {\rm O}(z^{2}) ,
\label{ops_Thexp:b}
\end{multline}
\begin{multline}
  \Theta^*_n(z)=
 -(n+\sum^m_{j=1}\rho_j)\frac{\bar{\phi}_{n}(0)}{\bar{\phi}_{n+1}(0)}z^{m-2} 
 \\
 + \bigg\{  [(n+\sum^m_{j=1}\rho_j)\sum^m_{j=1}z_j - \sum^m_{j=1}\rho_j z_j]
             \frac{\bar{\phi}_{n}(0)}{\bar{\phi}_{n+1}(0)}
            +(n+1+\sum^m_{j=1}\rho_j)\frac{\bar{\phi}_{n}(0)}{\bar{\phi}_{n+1}(0)}
             \frac{l_{n+1}}{\kappa_{n+1}}
 \\ 
 -(n-1+\sum^m_{j=1}\rho_j)\frac{\kappa_n\bar{\phi}_{n-1}(0)+\bar{\phi}_{n}(0)l_{n-1}
                               }{\bar{\phi}_{n+1}(0)}
      \bigg\}z^{m-3}
 \\ + {\rm O}(z^{m-4}) ,
\label{ops_ThSexp:a}
\end{multline}
\begin{multline}
  \Theta^*_n(z) =
 -[2V(0)-(n+1)W'(0)]\frac{\kappa_n}{\kappa_{n+1}}
 \\
 + \bigg\{-[2V'(0)-\shalf (n+1)W''(0)]\frac{\kappa_{n}}{\kappa_{n+1}}
           -[2V(0)-nW'(0)]\frac{\bar{l}_n}{\kappa_{n+1}}
 \\
 + [(n+2)W'(0)-2V(0)]\left( 
     \frac{\kappa^3_{n}}{\kappa_{n+2}\kappa^2_{n+1}}
     \frac{\bar{\phi}_{n+2}(0)}{\bar{\phi}_{n+1}(0)}
                   -\frac{\kappa_n}{\kappa_{n+1}}\frac{\bar{l}_{n+1}}{\kappa_{n+1}}
                       \right)
   \bigg\}z
 \\ + {\rm O}(z^{2}) ,
\label{ops_ThSexp:b}
\end{multline}
\begin{multline}
  \Omega_n(z) =
 (1+\half\sum^m_{j=1}\rho_j)z^{m-1}
 \\
 + \bigg\{ -\half(\sum^m_{j=1}\rho_j)(\sum^m_{j=1}z_j)
            +\half \sum^m_{j=1}\rho_j z_j 
            -\sum^m_{j=1}z_j
 \\
 + (n+2+\sum^m_{j=1}\rho_j)
           \frac{\kappa^2_n}{\kappa_{n+2}\kappa_{n+1}}
           \frac{\phi_{n+2}(0)}{\phi_{n+1}(0)}
     - \frac{l_{n+1}}{\kappa_{n+1}} \bigg\}z^{m-2}
 \\ + {\rm O}(z^{m-3}) ,
\label{ops_Omexp:a}
\end{multline}
\begin{multline}
  \Omega_n(z) =
 V(0)-nW'(0)
 \\
 + \bigg\{ V'(0)-\shalf nW''(0)
           +\left( V(0)\frac{\kappa_n}{\kappa_{n+1}}
                    +[V(0)-nW'(0)]\frac{\kappa_{n+1}}{\kappa_{n}} \right)
            \frac{\phi_{n}(0)}{\phi_{n+1}(0)}
 \\
 + [V(0)-nW'(0)]\frac{\bar{l}_n}{\kappa_n}
   - [V(0)-(n+1)W'(0)]\frac{\bar{l}_{n+1}}{\kappa_{n+1}} \bigg\}z
 \\ + {\rm O}(z^2) , 
\label{ops_Omexp:b}
\end{multline}
\begin{multline}
  \Omega^*_n(z) = 
 -\half\sum^m_{j=1}\rho_j z^{m-1}
 \\
 + \bigg\{  \half(\sum^m_{j=1}\rho_j)(\sum^m_{j=1}z_j)
            -\half \sum^m_{j=1}\rho_j z_j
            -(n+\sum^m_{j=1}\rho_j)\frac{\kappa_n}{\kappa_{n+1}}
             \frac{\bar{\phi}_{n}(0)}{\bar{\phi}_{n+1}(0)}
            +\frac{l_{n+1}}{\kappa_{n+1}} \bigg\}z^{m-2}
 \\
 + {\rm O}(z^{m-3}) ,
\label{ops_OmSexp:a}
\end{multline}
\begin{multline}
  \Omega^*_n(z) =
 (n+1)W'(0)-V(0)
 \\
 + \bigg\{ \shalf (n+1)W''(0)-V'(0) 
           +[(n+2)W'(0)-2V(0)]
            \frac{\kappa^2_{n}}{\kappa_{n+2}\kappa_{n+1}}
            \frac{\bar{\phi}_{n+2}(0)}{\bar{\phi}_{n+1}(0)}
           -W'(0)\frac{\bar{l}_{n+1}}{\kappa_{n+1}} \bigg\}z
 \\
 + {\rm O}(z^2) .
\label{ops_OmSexp:b}
\end{multline}
\end{proposition}

\begin{proof}
Following the approach of Laguerre \cite{La_1972b} we write $ F(z) $ in terms of
$ \phi_n(z) $, $ \psi_n(z) $ , $ \epsilon_n(z) $ and use (\ref{ops_FD}) to deduce
\begin{align}
  0 & = WF'-2VF-U ,
  \\
  & = W\left(\frac{\epsilon_n-\psi_n}{\phi_n}\right)' 
        -2V\frac{\epsilon_n-\psi_n}{\phi_n}-U ,
  \nonumber\\
  & = \frac{W(\psi_n\phi'_n-\phi_n\psi'_n)+2V\phi_n\psi_n-U\phi^2_n}{\phi^2_n}
       + W\left(\frac{\epsilon_n}{\phi_n}\right)'-2V\frac{\epsilon_n}{\phi_n} .
  \nonumber
\end{align}
The numerator of the first term is independent of $ \epsilon_n $, and so
is a polynomial in $ z $, and we denote this by
\begin{equation}
  2\frac{\phi_{n+1}(0)}{\kappa_n} z^n\Theta_n(z) =
  W(-\phi_n\epsilon'_n+\epsilon_n\phi'_n)+2V\phi_n\epsilon_n . 
\label{ops_Thdfn}
\end{equation}
Given that this is a polynomial we can determine its degree and minimum power 
of $ z $ by utilising the expansions of $ \phi_n, \epsilon_n $ both inside
and outside the unit circle, namely (\ref{ops_phiexp:a},\ref{ops_epsexp:a}).
We find the degree of the right-hand side is $ n+m-2 $ so that $ \Theta_n(z) $ 
is a polynomial of degree $ m-2 $. Developing the expansions further we arrive 
at (\ref{ops_Thexp:a}). An identical argument applies to the other combination
\begin{equation}
   2\frac{\bar{\phi}_{n+1}(0)}{\kappa_n} z^{n+1}\Theta^*_n(z) =
  W(\phi^*_n\epsilon^*_n{\!'}-\epsilon^*_n\phi^*_n{\!'})-2V\phi^*_n\epsilon^*_n ,
\label{ops_ThSdfn}
\end{equation}
and $ \Theta^*_n(z) $ is also a polynomial of degree $ m-2 $ with the expansion
(\ref{ops_Thexp:b}). To establish (\ref{ops_Omexp:a}) we utilise the other form 
of $ \Theta_n(z) $ and (\ref{ops_Cas:a}) to deduce
\begin{align}
  W(\psi_n\phi'_n-\phi_n\psi'_n)+2V\phi_n\psi_n-U\phi^2_n
  & = 2\frac{\phi_{n+1}(0)}{\kappa_n}z^n\Theta_n(z) ,
  \\
  & = [\phi_{n+1}\psi_n-\psi_{n+1}\phi_n]\Theta_n(z) .
  \nonumber
\end{align}
Separating those terms with $ \phi_n $ and $ \psi_n $ as factors we have
\begin{equation}
   \left\{ \Theta_n(z)\phi_{n+1}-W\phi'_n-V\phi_n \right\} \psi_n
 = \left\{ \Theta_n(z)\psi_{n+1}-W\psi'_n+V\psi_n-U\phi_n \right\} \phi_n ,
\end{equation}
so that this polynomial contains both $ \phi_n $ and $ \psi_n $ as factors and 
can be written as $ \Omega_n\phi_n\psi_n $ with $ \Omega_n(z) $ a polynomial 
of bounded degree. This latter polynomial can be defined as
\begin{align}
  2\frac{\phi_{n+1}(0)}{\kappa_n} z^n\Omega_n(z) 
  & = W(\psi_{n+1}\phi'_n-\phi_{n+1}\psi'_n)
      +V(\phi_n\psi_{n+1}+\psi_n\phi_{n+1})-U\phi_n\phi_{n+1} ,
  \nonumber\\
  & = W(\epsilon_{n+1}\phi'_n-\phi_{n+1}\epsilon'_n)
      +V(\phi_n\epsilon_{n+1}+\epsilon_n\phi_{n+1}) .
\label{ops_Omdfn}
\end{align}
Again employing the expansions (\ref{ops_phiexp:a},\ref{ops_epsexp:a}) we
determine the degree of $ \Omega_n(z) $ to be $ m-1 $ and the expansion
(\ref{ops_Omexp:a}) follows. Starting with the alternative definition of
$ \Theta^*_n(z) $ and (\ref{ops_Cas:b})
\begin{align}
  W(\phi^*_n\psi^*_n{\!'}-\psi^*_n\phi^*_n{\!'})-2V\phi^*_n\psi^*_n-U\phi^{*2}_n
  & = 2\frac{\bar{\phi}_{n+1}(0)}{\kappa_n} z^{n+1}\Theta^*_n(z) ,
  \\
  & = [\phi^*_{n+1}\psi^*_n-\psi^*_{n+1}\phi^*_n]\Theta^*_n(z) .
  \nonumber
\end{align}
and using the above argument we identify for the polynomial $ \Omega^*_n(z) $
\begin{align}
  2\frac{\bar{\phi}_{n+1}(0)}{\kappa_n} z^{n+1}\Omega^*_n(z) 
  & = W(-\psi^*_{n+1}\phi^*_n{\!'}+\phi^*_{n+1}\psi^*_n{\!'})
      -V(\phi^*_n\psi^*_{n+1}+\psi^*_n\phi^*_{n+1})
  \nonumber\\
  & \phantom{=}\qquad -U\phi^*_n\phi^*_{n+1} ,
  \nonumber\\
  & = W(-\epsilon^*_{n+1}\phi^*_n{\!'}+\phi^*_{n+1}\epsilon^*_n{\!'})
      -V(\phi^*_n\epsilon^*_{n+1}+\epsilon^*_n\phi^*_{n+1}) .
\label{ops_OmSdfn}
\end{align}
The degree of $ \Omega^*_n(z) $ is found to be $ m-1 $ and it has the expansion
(\ref{ops_OmSexp:a}). 
\end{proof}

\begin{remark}
Solving for $ \phi'_n $ and $ \epsilon'_n $ between (\ref{ops_Thdfn}) and 
(\ref{ops_Omdfn}) leads to (\ref{ops_zD:a}) and (\ref{ops_zD:b}),
whilst solving for $ \phi^*_n{\!'} $ and $ \epsilon^*_n{\!'} $ using
(\ref{ops_ThSdfn},\ref{ops_OmSdfn}) yields (\ref{ops_zD:c}) and 
(\ref{ops_zD:d}).
\end{remark}

Bilinear residue formulae relating products of a polynomial and an associated 
function evaluated at a finite singular point will arise in the theory of the 
deformation derivatives later. These are consequences of the workings of the 
proof of Proposition \ref{ops_SCpoly}, so we give a complete list presently.
  
\begin{corollary}
Bilinear residues are related to the coefficient function residues in the 
following equations, valid for all $ z_j $
\begin{align}
   \phi_n(z_j)\epsilon_n(z_j) 
  & =  2\frac{\phi_{n+1}(0)}{\kappa_n}z^n_j\frac{\Theta_n(z_j)}{2V(z_j)} ,
  \label{ops_BilRes:a} \\
   \phi^*_n(z_j)\epsilon^*_n(z_j) 
  & = -2\frac{\bar{\phi}_{n+1}(0)}{\kappa_n}z^{n+1}_j\frac{\Theta^*_n(z_j)}{2V(z_j)} ,
  \label{ops_BilRes:b}
\end{align}
\begin{align}
   \phi_{n+1}(z_j)\epsilon_n(z_j) 
  & =  2\frac{\phi_{n+1}(0)}{\kappa_n}z^n_j\frac{\Omega_n(z_j)+V(z_j)}{2V(z_j)} ,
  \label{ops_BilRes:c} \\
   \phi_{n}(z_j)\epsilon_{n+1}(z_j) 
  & =  2\frac{\phi_{n+1}(0)}{\kappa_n}z^n_j\frac{\Omega_n(z_j)-V(z_j)}{2V(z_j)} ,
  \label{ops_BilRes:d} \\
   \phi^*_{n}(z_j)\epsilon^*_{n+1}(z_j) 
  & = -2\frac{\bar{\phi}_{n+1}(0)}{\kappa_n}z^{n+1}_j
        \frac{\Omega^*_n(z_j)+V(z_j)}{2V(z_j)} ,
  \label{ops_BilRes:e} \\
   \phi^*_{n+1}(z_j)\epsilon^*_{n}(z_j) 
  & = -2\frac{\bar{\phi}_{n+1}(0)}{\kappa_n}z^{n+1}_j
        \frac{\Omega^*_n(z_j)-V(z_j)}{2V(z_j)} ,
  \label{ops_BilRes:f} \\
   \phi_n(z_j)\epsilon^*_n(z_j) 
  & = -\frac{z^n_j}{V(z_j)}\left[ \Omega_n(z_j)-V(z_j)
        -\frac{\kappa_{n+1}}{\kappa_n}z_j\Theta_n(z_j) \right] ,
  \label{ops_BilRes:g} \\
  & = -\frac{z^n_j}{V(z_j)}\left[ \Omega^*_n(z_j)-V(z_j)
        -\frac{\kappa_{n+1}}{\kappa_n}\Theta^*_n(z_j) \right] ,
  \label{ops_BilRes:h} \\
   \phi^*_n(z_j)\epsilon_n(z_j) 
  & =  \frac{z^n_j}{V(z_j)}\left[ \Omega_n(z_j)+V(z_j)
        -\frac{\kappa_{n+1}}{\kappa_n}z_j\Theta_n(z_j) \right] ,
  \label{ops_BilRes:i} \\
  & =  \frac{z^n_j}{V(z_j)}\left[ \Omega^*_n(z_j)+V(z_j)
        -\frac{\kappa_{n+1}}{\kappa_n}\Theta^*_n(z_j) \right] .
  \label{ops_BilRes:j}
\end{align}
\end{corollary}
\begin{proof}
These are all found by evaluating one of (\ref{ops_Thdfn}), (\ref{ops_ThSdfn}),
(\ref{ops_Omdfn}), or (\ref{ops_OmSdfn}) at $ z=z_j $ and using 
(\ref{ops_Cas:a}-\ref{ops_Cas:c}).
\end{proof}

The initial members of the sequences of coefficient functions 
$ \{\Theta_n\}^{\infty}_{n=0} $, $ \{\Theta^*_n\}^{\infty}_{n=0} $,
$ \{\Omega_n\}^{\infty}_{n=0} $, $ \{\Omega^*_n\}^{\infty}_{n=0} $ are given by
\begin{align}
  2\frac{\phi_1(0)}{\kappa_0}\Theta_0(z) = & 2V(z)-\kappa^2_0 U(z) ,
  \label{ops_Theta0} \\
  2\frac{\phi_2(0)}{\kappa_1}z\Theta_1(z) = 
   & \frac{\kappa^2_1}{\kappa^2_0} z^2(2V(z)-\kappa^2_0 U(z)) 
     - 2\kappa_1\phi_1(0)zU(z)
  \label{ops_Theta1} \\
   & \qquad -2\frac{\kappa_1\phi_1(0)}{\kappa^2_0}W(z) 
     - \frac{\phi^2_1(0)}{\kappa^2_0}(2V(z)+\kappa^2_0 U(z)) ,
  \nonumber \\
  2\frac{\bar{\phi}_1(0)}{\kappa_0}z\Theta^*_0(z) = & -2V(z)-\kappa^2_0 U(z) ,
  \label{ops_ThetaS0} \\
  2\frac{\bar{\phi}_2(0)}{\kappa_1}z^2\Theta^*_1(z) = 
   & \frac{\bar{\phi}^2_1(0)}{\kappa^2_0} z^2(2V(z)-\kappa^2_0 U(z)) 
     - 2\kappa_1\bar{\phi}_1(0)zU(z)
  \label{ops_ThetaS1} \\
   & \qquad -2\frac{\kappa_1\bar{\phi}_1(0)}{\kappa^2_0}W(z) 
     - \frac{\kappa^2_1}{\kappa^2_0} (2V(z)+\kappa^2_0 U(z)) ,
  \nonumber
\end{align}
\begin{align}
  2\phi_1(0)\Omega_0(z) = & \kappa_1z(2V(z)-\kappa^2_0 U(z)) - \kappa^2_0\phi_1(0)U(z) ,
  \label{ops_Omega0} \\
  2\bar{\phi}_1(0)z\Omega^*_0(z) = & -\kappa_1(2V(z)+\kappa^2_0 U(z)) 
     - \kappa^2_0\bar{\phi}_1(0)zU(z) .
  \label{ops_OmegaS0}
\end{align}
 
One can take combinations of the above functional-difference equations and 
construct exact differences when $ z $ is evaluated at the finite singular points 
of the weight, i.e. $ W(z)=0 $. The integration of the system is given in the
following proposition.

\begin{proposition}\label{ops_Bilinear}
At all the finite singular points $ z_j, j=1,\ldots m $, with the exception of $ z_j=0 $,
the coefficient functions satisfy the bilinear identities
\begin{gather}
  \Omega^2_n(z_j)
  = \frac{\kappa_n \phi_{n+2}(0)}{\kappa_{n+1} \phi_{n+1}(0)}z_j
      \Theta_n(z_j)\Theta_{n+1}(z_j)+V^2(z_j) ,
  \label{ops_OTeq:a} \\
  \Omega^{*2}_n(z_j)
  = \frac{\kappa_n \bar{\phi}_{n+2}(0)}{\kappa_{n+1} \bar{\phi}_{n+1}(0)}z_j
      \Theta^*_n(z_j)\Theta^*_{n+1}(z_j)+V^2(z_j) ,
  \label{ops_OTeq:b} \\
  \left[\Omega_{n-1}(z_j)
    -\frac{\kappa^2_{n-1}}{\kappa^2_{n}}\frac{\phi_{n+1}(0)}{\phi_{n}(0)}\Theta_n(z_j)
  \right]^2
  = \frac{\phi_{n+1}(0)\bar{\phi}_{n}(0)}{\kappa^2_{n}}
      \Theta_n(z_j)\Theta^*_{n-1}(z_j)+V^2(z_j) ,
  \label{ops_OTeq:c} \\
  \left[\Omega^*_{n-1}(z_j)
    -\frac{\kappa^2_{n-1}}{\kappa^2_{n}}\frac{\bar{\phi}_{n+1}(0)}{\bar{\phi}_{n}(0)}
      z_j\Theta^*_n(z_j) ,
  \right]^2 \nonumber\\
  = \frac{\kappa_{n-1}\bar{\phi}_{n+1}(0)\phi_{n}(0)}{\kappa^3_{n}}
      z^2_j\Theta^*_n(z_j)\Theta_{n-1}(z_j)+V^2(z_j) ,
  \label{ops_OTeq:d} \\
    \frac{\phi_{n+1}(0)\bar{\phi}_{n+1}(0)}{\kappa^2_{n}}
      z_j\Theta_n(z_j)\Theta^*_{n}(z_j)+V^2(z_j)
  = \left[\Omega_{n}(z_j)-\frac{\kappa_{n+1}}{\kappa_{n}}z_j\Theta_n(z_j)\right]^2 ,
  \label{ops_OTeq:e} \\
  \phantom{
    \frac{\phi_{n+1}(0)\bar{\phi}_{n+1}(0)}{\kappa^2_{n}}
      z_j\Theta_n(z_j)\Theta^*_{n}(z_j)+V^2(z_j) }
  = \left[\Omega^*_{n}(z_j)-\frac{\kappa_{n+1}}{\kappa_{n}}\Theta^*_n(z_j)\right]^2 .
  \label{ops_OTeq:f}
\end{gather}
\end{proposition}
\begin{proof}[First Proof]
We take the first pair of identities (\ref{ops_OTeq:a}) and (\ref{ops_OTeq:b})
as an example for our first proof.
Multiplying the $ \Omega_n, \Omega_{n-1} $ terms of (\ref{ops_rrCf:a}) by the
corresponding terms of (\ref{ops_rrCf:b}), evaluated at a finite singular point
$ z=z_j $, one has an exact difference
\begin{multline}
  \Omega^2_n(z_j)-\Omega^2_{n-1}(z_j) \\
 =\frac{\kappa_n \phi_{n+2}(0)}{\kappa_{n+1} \phi_{n+1}(0)}z_j
      \Theta_n(z_j)\Theta_{n+1}(z_j)
 -\frac{\kappa_{n-1} \phi_{n+1}(0)}{\kappa_n \phi_{n}(0)}z_j
      \Theta_{n-1}(z_j)\Theta_{n}(z_j) ,
\end{multline}
assuming none of the $ z_j $ coincide with $ -r_{n+1}/r_n $ for any $ n $.
Upon summing this relation the summation constant is calculated to be
\begin{equation}
   \Omega^2_0(z_j)-\frac{\kappa_0 \phi_{2}(0)}{\kappa_{1} \phi_{1}(0)}z_j
      \Theta_0(z_j)\Theta_{1}(z_j)
   = V^2(z_j) ,
\end{equation}
by using the initial members of the coefficient function sequences in 
(\ref{ops_Omega0},\ref{ops_Theta0},\ref{ops_Theta1}). The result is 
(\ref{ops_OTeq:a}), whilst the second relation follows from an identical argument 
applied to (\ref{ops_rrCf:c},\ref{ops_rrCf:d}).
\end{proof}

\begin{proof}[Second Proof]
The three pairs of formulae (\ref{ops_OTeq:a},\ref{ops_OTeq:b}), 
(\ref{ops_OTeq:c},\ref{ops_OTeq:d}) and (\ref{ops_OTeq:e},\ref{ops_OTeq:f})
arise from the fact that at a finite singular point $ z_j $ the determinant of the 
matrix spectral derivative must vanish. Thus (\ref{ops_OTeq:a}) and 
(\ref{ops_OTeq:b}) express the condition that the determinant of the matrix on
the right-hand sides of (\ref{ops_XzDer:a}) and (\ref{ops_XzDer:b}) vanish
respectively. It can be shown that the same condition applied to the right-hand 
sides of (\ref{ops_ZzDer:b}) and (\ref{ops_ZzDer:a}) implies (\ref{ops_OTeq:c}) 
and (\ref{ops_OTeq:d}) respectively when one takes into account the identities 
(\ref{ops_rrCf:f}), (\ref{ops_rrCf:i}), (\ref{ops_rrCf:j}) and (\ref{ops_rrCf:c}).
The last pair are a consequence of $ \det(WA_n(z_j;t)) = 0 $ along with the 
identity (\ref{ops_rrCf:j}). 
\end{proof}

\begin{proof}[Third Proof]
All the bilinear identities in Proposition \ref{ops_Bilinear} can be easily
derived from the residue formulae (\ref{ops_BilRes:a}-\ref{ops_BilRes:j}) by 
multiplying any two of the above formulae and then factoring the resulting 
product in a different way. Thus (\ref{ops_OTeq:a}) arises from multiplying
(\ref{ops_BilRes:c}) and (\ref{ops_BilRes:d}) and then factoring the product
in order to employ (\ref{ops_BilRes:a}). Equation (\ref{ops_OTeq:c}) comes from
multiplying (\ref{ops_BilRes:a}) and (\ref{ops_BilRes:b}) with $ n \mapsto n-1 $,
using the recurrences (\ref{ops_rr:b}), (\ref{ops_rre:b}) with $ n \mapsto n-1 $
to solve for $ \phi^*_{n-1}(z_j), \epsilon^*_{n-1}(z_j) $ and employing 
(\ref{ops_BilRes:a}) along with (\ref{ops_BilRes:c}) and (\ref{ops_BilRes:d}) setting 
$ n \mapsto n-1 $. Equation (\ref{ops_OTeq:e}) is derived by multiplying 
(\ref{ops_BilRes:a}) and (\ref{ops_BilRes:b}) and then factoring using
(\ref{ops_BilRes:g}) and (\ref{ops_BilRes:i}). The reciprocal versions follow
from similar reasoning.
\end{proof}

\begin{remark}
It is clear from the first proof that the bilinear identities given in Proposition 
\ref{ops_Bilinear} can be straightforwardly generalised to ones that are functions 
of $ z $ rather than evaluated at special $ z$ values. 
They can be derived directly from Proposition \ref{ops_Linear1}, so apply in situations where
the weights are not semi-classical, and contain additional
terms with a factor of $ W(z) $ and sums of products of other coefficients ranging from 
$ j = 1,\ldots n $. 
However because we will have no use for such relations we refrain from writing these down.
\end{remark}

\begin{remark}
As $ z=0 $ is a finite singular point then the limit as $ z \to 0 $ may be taken in the
product of (\ref{ops_rrCf:a},\ref{ops_rrCf:b}), however this does not lead to
any new independent relation but simply recovers
\begin{equation*}
   \Omega_n(0) = V(0)-nW'(0) .
\end{equation*}
\end{remark}

In the case of a regular semi-classical weight function the matrix $ A_n(z;t) $ 
has the partial fraction decomposition
\begin{equation}
   A_n(z;t) := \sum^{m}_{j=1}\frac{A_{nj}}{z-z_j} ,
\label{An_parfrac}
\end{equation}
under the assumptions following (\ref{ops_scwgt}). Let us take the first finite 
singularity to be situated at the origin, $ z_1=0 $. Some care needs to exercised 
as many relations differ depending on whether $ z_j=0 $ or not because of the 
additional term in (\ref{ops_rrCf:j}).
The residue matrices for the finite singularities are given by
\begin{multline}
   A_{nj} = \frac{\rho_j}{2V(z_j)} \\
   \times
       \begin{pmatrix}
              -\Omega_n(z_j)-V(z_j)
              +\dfrac{\kappa_{n+1}}{\kappa_n}z_j\Theta_n(z_j)
            & \dfrac{\phi_{n+1}(0)}{\kappa_n}\Theta_n(z_j)
            \cr
              -\dfrac{\bar{\phi}_{n+1}(0)}{\kappa_n}z_j\Theta^*_n(z_j)
            &  \Omega^*_n(z_j)-V(z_j)
                     -\dfrac{\kappa_{n+1}}{\kappa_n}\Theta^*_n(z_j)
            \cr
       \end{pmatrix} ,
\label{An_resJ}
\end{multline}
for $ j=2,\ldots,m $ and 
\begin{equation}
   A_{n1} = \frac{\rho_1}{2V(0)}
       \begin{pmatrix}
              nW'(0)-2V(0) & [2V(0)-nW'(0)]r_n \cr 0 & 0 \cr
       \end{pmatrix} .
\label{An_res0}
\end{equation}
Using the identity (\ref{ops_rrCf:j}) we note that 
\begin{align}
  {\rm Tr} A_{nj} & = -\rho_j, \quad j=2,\ldots,m, 
  \label{An_trJ} \\
  {\rm Tr} A_{n1} & = n-\rho_1 . 
  \label{An_tr0}
\end{align}
In either case we find that $ \det A_{nj} = 0 $ using (\ref{ops_OTeq:e}).
An alternative expression for the residue matrices in the case $ j=2,\ldots,m $
is
\begin{equation}
   A_{nj} = -\half\rho_j z_j^{-n}
       \begin{pmatrix}
              \phi^*_n(z_j)\epsilon_n(z_j) & -\phi_n(z_j)\epsilon_n(z_j) \cr
              -\phi^*_n(z_j)\epsilon^*_n(z_j) & \phi_n(z_j)\epsilon^*_n(z_j) \cr
       \end{pmatrix} .
\label{An_resBil}
\end{equation}
The regular singularity at $ z=\infty $ has a residue matrix given by
\begin{equation}
  A_{n\infty} := {\rm Res}_{z=0}\left\{ -z^{-2}A_n(z^{-1}) \right\}
               = -\sum^{m}_{j=1} A_{nj} .  
\label{An_inftyDef}
\end{equation}
Using this definition and the large $ z $ terms for the coefficient functions
(\ref{ops_Thexp:a}), (\ref{ops_ThSexp:a}), (\ref{ops_Omexp:a}), (\ref{ops_OmSexp:a})
we evaluate this matrix to be
\begin{equation}
   A_{n\infty} = 
       \begin{pmatrix}
              -n & 0 \cr
              -(n+\sum^{m}_{j=1}\rho_j)\bar{r}_n & \sum^{m}_{j=1}\rho_j \cr
       \end{pmatrix} .
\label{An_infty}
\end{equation}
We read off that $ {\rm Tr} A_{n\infty} = -n+\sum^{m}_{j=1}\rho_j $ and 
$ \det A_{n\infty} = -n\sum^{m}_{j=1}\rho_j $. Furthermore relations (\ref{An_infty}) and
(\ref{An_inftyDef}) imply the summation identities
\begin{align} 
   \half\sum^{m}_{j=1}\rho_jz_j^{-n}\phi_n(z_j)\epsilon_n(z_j)
   & = 0 ,\\
   \half\sum^{m}_{j=1}\rho_jz_j^{-n}\phi^*_n(z_j)\epsilon_n(z_j)
   & = -n ,\\
   \half\sum^{m}_{j=1}\rho_jz_j^{-n}\phi_n(z_j)\epsilon^*_n(z_j)
   & = \sum^m_{j=1}\rho_j ,\\
   \half\sum^{m}_{j=1}\rho_jz_j^{-n}\phi^*_n(z_j)\epsilon^*_n(z_j)
   & = (n+\sum^m_{j=1}\rho_j)\bar{r}_n .
\end{align} 

We wish to close this part by commenting on how one would obtain the discrete analogs
of the Schlesinger equations or multi-variable extensions of the discrete Painlev\'e
equations from the theory outlined above. 
If one evaluates (\ref{ops_OTeq:a}) (or (\ref{ops_OTeq:b}) for that matter)
at two distinct singularities $ z_1, z_2 $, consolidates terms and then takes their
ratio the result is
\begin{equation*}
   \frac{z_1\Theta_n(z_1)\Theta_{n+1}(z_1)}{z_2\Theta_n(z_2)\Theta_{n+1}(z_2)}
   = \frac{[\Omega_n(z_1)-V(z_1)][\Omega_n(z_1)+V(z_1)]}
          {[\Omega_n(z_2)-V(z_2)][\Omega_n(z_2)+V(z_2)]} .
\end{equation*} 
This constitutes a recurrence relation for $ \Theta_n(z_1)/\Theta_n(z_2) $.
To find a recurrence involving $ \Omega_n $ one adopts another method. 
By comparing the expansions inside and outside the unit circle, for example
equations (\ref{ops_Thexp:a}) and (\ref{ops_Thexp:b}) for $ \Theta_n(z) $, and in
particular where they overlap one can derive expressions for the sub-leading 
coefficients $ l_n, m_n $ in terms of the higher ones. Employing these expressions
in one of the expansion forms for $ \Omega_n(z_{1,2}) $ ((\ref{ops_Omexp:a}) or 
(\ref{ops_Omexp:b})) will yield a recurrence for $ \Omega_n $.

We now consider the dynamics of deforming the semi-classical weight 
(\ref{ops_scwgt2}) through a $t$-dependence of the finite singular points $ z_j(t) $,
\begin{equation}
   \frac{\dot{w}}{w} = -\sum^m_{j=1}\rho_j\frac{\dot{z}_j}{z-z_j} ,
\label{ops_tDwgt}
\end{equation}
where $ \;\dot{} := d/dt $.
Given this motion of the finite singularities we consider the $t$-derivatives of
the bi-orthogonal polynomial system.
  
\begin{proposition}\label{ops_deformD}
The deformation derivative of a semi-classical bi-orthogonal polynomial is
\begin{multline}
  \dot{\phi}_n(z) =
  \Big\{-\frac{\dot{\kappa}_n}{\kappa_n}
        -\sum^m_{j=1}\rho_j\frac{\dot{z}_j}{z_j}
        +\half\sum^m_{j=1}\rho_j\frac{\dot{z}_j}{z_j}
              z^{-n}_j \epsilon_n(z_j)\phi^*_n(z_j)\frac{z}{z-z_j}
  \Big\}\phi_n(z) 
  \\
  -\Big\{\half\sum^m_{j=1}\rho_j\frac{\dot{z}_j}{z_j}
              z^{1-n}_j \epsilon_n(z_j)\phi_n(z_j)\frac{1}{z-z_j}
  \Big\}\phi^*_n(z) ,
\label{ops_tD:a}
\end{multline}
whilst that of its reciprocal polynomial is
\begin{multline}
  \dot{\phi}^*_n(z) =
  \Big\{-\frac{\dot{\kappa}_n}{\kappa_n}
        +\half\sum^m_{j=1}\rho_j\frac{\dot{z}_j}{z_j}
              z^{1-n}_j \epsilon^*_n(z_j)\phi_n(z_j)\frac{1}{z-z_j}
  \Big\}\phi^*_n(z) 
  \\
  -\Big\{\half\sum^m_{j=1}\rho_j\frac{\dot{z}_j}{z_j}
              z^{-n}_j \epsilon^*_n(z_j)\phi^*_n(z_j)\frac{z}{z-z_j}
  \Big\}\phi_n(z) .
\label{ops_tD:b}
\end{multline}
The deformation derivative of an associated function is
\begin{multline}
  \dot{\epsilon}_n(z) =
  \Big\{-\frac{\dot{\kappa}_n}{\kappa_n}
        -\half\sum^m_{j=1}\rho_j\frac{\dot{z}_j}{z_j}
              z^{-n}_j \epsilon^*_n(z_j)\phi_n(z_j)\frac{z}{z-z_j}
  \Big\}\epsilon_n(z) 
  \\
  +\Big\{\half\sum^m_{j=1}\rho_j\frac{\dot{z}_j}{z_j}
              z^{1-n}_j \epsilon_n(z_j)\phi_n(z_j)\frac{1}{z-z_j}
  \Big\}\epsilon^*_n(z) ,
\label{ops_tD:c}
\end{multline}
and that of a reciprocal associated function is
\begin{multline}
  \dot{\epsilon}^*_n(z) =
  \Big\{-\frac{\dot{\kappa}_n}{\kappa_n}
        -\half\sum^m_{j=1}\rho_j\frac{\dot{z}_j}{z_j}
              z^{1-n}_j \epsilon_n(z_j)\phi^*_n(z_j)\frac{1}{z-z_j}
  \Big\}\epsilon^*_n(z) 
  \\
  +\Big\{\half\sum^m_{j=1}\rho_j\frac{\dot{z}_j}{z_j}
              z^{-n}_j \epsilon^*_n(z_j)\phi^*_n(z_j)\frac{z}{z-z_j}
  \Big\}\epsilon_n(z) .
\label{ops_tD:d}
\end{multline}
\end{proposition}

\begin{proof}
Differentiating the orthonormality condition 
\begin{equation*}
   \int\frac{d\zeta}{2\pi i\zeta}w(\zeta)
               \phi_{n}(\zeta)\bar{\phi}_{n-i}(\bar{\zeta}) = \delta_{i,0} ,
\end{equation*}
and using (\ref{ops_tDwgt}) we find
\begin{equation*}
  0 = \frac{\dot{\kappa}_n}{\kappa_n}\delta_{i,0}
      + \int\frac{d\zeta}{2\pi i\zeta}w(\zeta)\dot{\phi}_n\bar{\phi}_{n-i}
      - \sum_j\rho_j \dot{z}_j \int\frac{d\zeta}{2\pi i\zeta}w(\zeta)
        \frac{1}{\zeta-z_j}\phi_n\bar{\phi}_{n-i}, \; i=0,\ldots,n .
\end{equation*}
Now 
\begin{align*}
  \int\frac{d\zeta}{2\pi i\zeta}w(\zeta)\frac{1}{\zeta-z}
  \phi_n(\zeta)\bar{\phi}_{n-i}(\bar{\zeta})
  & = \int\frac{d\zeta}{2\pi i\zeta}w(\zeta)\phi_n(\zeta)
          \frac{\bar{\phi}_{n-i}(\zeta^{-1})-\bar{\phi}_{n-i}(z^{-1})}{\zeta-z}
  \\
  & \qquad
      + \bar{\phi}_{n-i}(z^{-1})
      \int\frac{d\zeta}{2\pi i\zeta}w(\zeta)\frac{\phi_{n}(\zeta)}{\zeta-z} ,
  \\
  & = -\frac{1}{z}\delta_{i,0}
      + \frac{\bar{\phi}_{n-i}(z^{-1})}{2z}\epsilon_n(z), \quad n > 0 ,
\end{align*}
so that 
\begin{equation*}
  0 = \left(\frac{\dot{\kappa}_n}{\kappa_n}+\sum_j \rho_j\frac{\dot{z}_j}{z_j}
      \right)\delta_{i,0}
  - \half\sum^m_{j=1}\rho_j\frac{\dot{z}_j}{z_j}
    z^{i-n}_j\phi^*_{n-i}(z_j)\epsilon_n(z_j)
  + \int\frac{d\zeta}{2\pi i\zeta}w(\zeta)\dot{\phi}_n\bar{\phi}_{n-i} .
\end{equation*}
In addition we can represent $ \bar{\phi}_{n-i}(z) $ as
\begin{align*}
  \bar{\phi}_{n-i}(z) & = \sum^n_{j=0} \delta_{i,j}\bar{\phi}_{n-j}(z) ,
  \\
  & = \sum^n_{j=0} \int\frac{d\zeta}{2\pi i\zeta}w(\zeta)
      \bar{\phi}_{n-i}(\bar{\zeta)}\phi_{n-j}(\zeta)\bar{\phi}_{n-j}(z) ,
  \\
  &  = \int\frac{d\zeta}{2\pi i\zeta}w(\zeta)\bar{\phi}_{n-i}(\bar{\zeta})
       \sum^n_{j=0} \phi_{n-j}(\zeta)\bar{\phi}_{n-j}(z) ,
  \\
  & = \int\frac{d\zeta}{2\pi i\zeta}w(\zeta)\bar{\phi}_{n-i}(\bar{\zeta})
      \frac{\phi^*_n(\zeta)\overline{\phi^*_n}(z)
       -\zeta z\phi_n(\zeta)\bar{\phi}_n(z)}{1-\zeta z} .
\end{align*}
Writing the Kronecker delta in a similar way the whole expression becomes
\begin{multline*}
  0 = \int\frac{d\zeta}{2\pi i\zeta}w(\zeta)\bar{\phi}_{n-i}(\bar{\zeta})
      \bigg\{ \dot{\phi}_n(\zeta) 
            +\left(\frac{\dot{\kappa}_n}{\kappa_n}
                   +\sum^m_{j=1}\rho_j\frac{\dot{z}_j}{z_j}\right)\phi_n(\zeta)
  \\
  -\half\sum^m_{j=1}\rho_j\frac{\dot{z}_j}{z_j}\epsilon_n(z_j)
   \frac{\phi^*_n(\zeta)z^{-n}_j\phi_n(z_j)
       -\zeta z^{-1}_j\phi_n(\zeta)z^{-n}_j\phi^*_n(z_j)}{1-\zeta z^{-1}_j}
      \bigg\} ,
\end{multline*}
for all $ 0 \leq i \leq n $ and (\ref{ops_tD:a}) then follows. The second 
relation follows by an identical argument applied to 
\begin{equation*}
   \int\frac{d\zeta}{2\pi i\zeta}w(\zeta)
               \phi_{n-i}(\zeta)\bar{\phi}_{n}(\bar{\zeta}) = \delta_{i,0} .
\end{equation*}
The derivatives of the associated functions (\ref{ops_tD:c}), (\ref{ops_tD:d})
follow from differentiating the definitions (\ref{ops_eps:a}), (\ref{ops_eps:b})
and employing the first two results of the proposition along with the relation
(\ref{ops_Cas:c}).
\end{proof}

\begin{corollary}\label{cor_rdot}
The $t$-derivatives of the reflection coefficients are 
\begin{align}
   \frac{\dot{r}_n}{r_n} & =
   \half\sum^m_{j=1}\rho_j\frac{\dot{z}_j}{z_j}
   \frac{\Omega_{n-1}(z_j)-V(z_j)}{V(z_j)} , 
   \label{ops_rdot} \\
   \frac{\dot{\bar{r}}_n}{\bar{r}_n} & =
   \half\sum^m_{j=1}\rho_j\frac{\dot{z}_j}{z_j}
   \frac{\Omega^*_{n-1}(z_j)+V(z_j)}{V(z_j)} . 
   \label{ops_rCdot}
\end{align}
\end{corollary}

\begin{proof}
An alternative formula to (\ref{ops_tD:a}) is
\begin{equation}
  \dot{\phi}_n(\zeta) =
  -\left(\frac{\dot{\kappa}_n}{\kappa_n}
         +\sum^m_{j=1}\rho_j\frac{\dot{z}_j}{z_j}\right)\phi_n(\zeta)
  +\half\sum^m_{j=1}\rho_j\frac{\dot{z}_j}{z_j}\epsilon_n(z_j)
   \sum^n_{l=0} \bar{\phi}_{n-l}(z^{-1}_j)\phi_{n-l}(\zeta) ,
\end{equation}
and by examining the coefficients of $ \zeta^n, \zeta^0 $ we deduce that
\begin{equation*}
   \frac{\dot{r}_n}{r_n} = \half\frac{\kappa_{n-1}}{\phi_n(0)}
    \sum^m_{j=1}\rho_j\frac{\dot{z}_j}{z_j}
    z^{1-n}_j\epsilon_n(z_j)\phi_{n-1}(z_j) .
\end{equation*}
Noting that the derivative term of (\ref{ops_Omdfn}) vanishes when
$ z=z_j $ and employing (\ref{ops_Cas:a}) we arrive at (\ref{ops_rdot}).
The second equation, (\ref{ops_rCdot}), follows by identical reasoning.
\end{proof}

Sums of the bilinear residues over the finite singular points are related to deformation 
derivatives in the following way,
\begin{align}
  2\frac{\dot{\kappa}_n}{\kappa_n} 
  & = -\sum^m_{j=1}\rho_j\frac{\dot{z}_j}{z_j}
    + \half\sum^m_{j=1}\rho_j\frac{\dot{z}_j}{z_j}
                       z^{-n}_j\epsilon_n(z_j)\phi^*_{n}(z_j) ,
  \nonumber \\
  & = - \half\sum^m_{j=1}\rho_j\frac{\dot{z}_j}{z_j}
                       z^{-n}_j\epsilon^*_n(z_j)\phi_{n}(z_j) ,
  \label{ops_tDkappa} \\
  \frac{\dot{\phi}_n(0)}{\phi_n(0)}+\frac{\dot{\kappa}_n}{\kappa_n}
  + \sum^m_{j=1}\rho_j\frac{\dot{z}_j}{z_j}
  & = \half\frac{\kappa_n}{\phi_n(0)}\sum^m_{j=1}\rho_j\frac{\dot{z}_j}{z_j}
                       z^{-n}_j\epsilon_n(z_j)\phi_{n}(z_j) ,
  \nonumber \\
  & = \frac{\phi_{n+1}(0)}{\phi_n(0)}
      \sum^m_{j=1}\frac{\rho_j}{2V(z_j)}\frac{\dot{z}_j}{z_j}\Theta_n(z_j) ,
  \label{ops_tDphi} \\
  \frac{\dot{\bar{\phi}}_n(0)}{\bar{\phi}_n(0)}+\frac{\dot{\kappa}_n}{\kappa_n}
  & = - \half\frac{\kappa_n}{\bar{\phi}_n(0)} \sum^m_{j=1}\rho_j\frac{\dot{z}_j}{z_j}
                       z^{-n}_j\epsilon^*_n(z_j)\phi^*_{n}(z_j) ,
  \nonumber \\
  & = \frac{\bar{\phi}_{n+1}(0)}{\bar{\phi}_n(0)}
      \sum^m_{j=1}\frac{\rho_j}{2V(z_j)}\dot{z}_j\Theta^*_n(z_j) .
  \label{ops_tDCphi}
\end{align}

For the regular semi-classical weights we can also formulate the system of
deformation derivatives as a $2\times 2$ matrix differential equation and
demonstrate that the system preserves the monodromy data with respect 
to the motion of the finite singularities $ z_j(t) $. 

\begin{corollary}\label{ops_DefDer}
The deformation derivatives for a system of regular semi-classical bi-orthogonal 
polynomials and associated functions (\ref{ops_tD:a}-\ref{ops_tD:d}) are 
equivalent to the matrix differential equation
\begin{equation}
   \dot{Y}_{n} := B_n Y_{n}
   = \left\{ B_{\infty} - \sum^{m}_{j=1}\frac{\dot{z}_j}{z-z_j}A_{nj}
        \right\} Y_{n} . 
\label{ops_YtDer}
\end{equation} 
where
\begin{equation}
  B_{\infty} = 
       \begin{pmatrix}
               \dfrac{\dot{\kappa}_{n}}{\kappa_n}
            & 0
            \cr
               \dfrac{\kappa_n\dot{\bar{\phi}}_{n}(0)+\dot{\kappa}_{n}\bar{\phi}_n(0)}{
               \kappa^2_n}
            & -\dfrac{\dot{\kappa}_{n}}{\kappa_n}
            \cr
       \end{pmatrix} .
\end{equation}
\end{corollary}
\begin{proof}
This follows from a partial fraction decomposition of the system 
(\ref{ops_tD:a}-\ref{ops_tD:d}) and using (\ref{ops_tDkappa},\ref{ops_tDCphi}). 
\end{proof}

\begin{remark}
For a special class of irregular semi-classical weights a result analogous to Corollary
\ref{ops_DefDer} has been given by Bertola, Eynard and Harnad \cite{BEH_2003b}.
\end{remark}

In the case of the pair (\ref{ops_Yrecur}), (\ref{ops_YtDer}) compatibility 
implies the relation
\begin{equation}
   \dot{K}_n = B_{n+1}K_n-K_nB_n ,
\end{equation}
however there are no new identities arising from this condition. Taking the 
$ (1,1) $-component of both sides of this equation we see that it is identically satisfied 
through the use of (\ref{ops_rrCf:f}) and (\ref{ops_tDCphi}). Or if we take the
$ (1,2) $-components then they are equal when use of made of (\ref{ops_rrCf:a}) and
(\ref{ops_tDkappa},\ref{ops_rdot}). In a similar way we find both sides of the 
$ (2,1) $-components are identical when we employ (\ref{ops_rrCf:c}) and 
(\ref{ops_tDCphi}). Finally the $ (2,2) $-components on both sides are the same 
after taking into account (\ref{ops_rrCf:h}) and (\ref{ops_tDkappa},\ref{ops_tDCphi}).

For the pair of linear differential relations (\ref{ops_YzDer}), (\ref{ops_YtDer}) 
compatibility leads us to the Schlesinger equations 
\begin{gather}
   \dot{A}_{nj} = \left[ B_{\infty},A_{nj} \right] 
   + \sum_{k \neq j}\frac{\dot{z}_j-\dot{z}_k}{z_j-z_k} \left[ A_{nk}, A_{nj} \right] ,
   \\
   \dot{A}_{n\infty} = \left[ B_{\infty},A_{n\infty} \right] .
\label{Schlesinger}
\end{gather}
Again there is not anything essentially new here, that couldn't be derived from the
system of deformation derivatives (\ref{ops_tD:a}-\ref{ops_tD:d}), but it is an
efficient way to compute the deformation derivatives of bilinear products.
Employing the explicit representations of our matrices $ A_{nj} $ we find the 
following independent derivatives in component form
\begin{multline}
  \frac{d}{dt}\frac{\rho_j}{2V(z_j)}
  \left[ \Omega_n(z_j)+V(z_j)-\frac{\kappa_{n+1}}{\kappa_{n}}z_j\Theta_n(z_j) \right] \\
  = -\frac{\rho_j}{2V(z_j)}\frac{\phi_{n+1}(0)}{\kappa^3_n}
     \frac{d}{dt}(\kappa_n\bar{\phi}_n(0))\Theta_n(z_j)
  \\
    - \frac{\rho_j}{2V(z_j)}\frac{\phi_{n+1}(0)\bar{\phi}_{n+1}(0)}{\kappa^2_{n}}
      \sum_{k \neq j}\frac{\rho_k}{2V(z_k)}\frac{\dot{z}_j-\dot{z}_k}{z_j-z_k}
      \left[ z_k\Theta^*_n(z_k)\Theta_n(z_j)-z_j\Theta_n(z_k)\Theta^*_n(z_j) \right] ,
  \label{ops_Schl:a}
\end{multline}
\begin{multline}
  \frac{d}{dt}\frac{\rho_j}{2V(z_j)}
  \frac{\phi_{n+1}(0)}{\kappa_{n}}\Theta_n(z_j)
  = \frac{\rho_j}{V(z_j)}\frac{\phi_{n+1}(0)}{\kappa_n}\bigg\{
     \frac{\dot{\kappa}_n}{\kappa_n}\Theta_n(z_j)
    + \sum_{k \neq j}\frac{\rho_k}{2V(z_k)}\frac{\dot{z}_j-\dot{z}_k}{z_j-z_k}
  \\ \times
      \left[ \Theta_n(z_k)
             \big[\Omega_n(z_j)-\frac{\kappa_{n+1}}{\kappa_{n}}z_j\Theta_n(z_j)\big]
            -\Theta_n(z_j)
             \big[\Omega_n(z_k)-\frac{\kappa_{n+1}}{\kappa_{n}}z_k\Theta_n(z_k)\big] \right]
    \bigg\} ,
  \label{ops_Schl:b}
\end{multline}
\begin{multline}
  \frac{d}{dt}\frac{\rho_j}{2V(z_j)}
  \frac{\bar{\phi}_{n+1}(0)}{\kappa_{n}}z_j\Theta^*_n(z_j)
  = \frac{\rho_j}{V(z_j)}\frac{\bar{\phi}_{n+1}(0)}{\kappa_n}\bigg\{
     - \frac{\dot{\kappa}_n}{\kappa_n}z_j\Theta^*_n(z_j)
  \\
     + \frac{1}{\kappa_n\bar{\phi}_{n+1}(0)}\frac{d}{dt}(\kappa_n\bar{\phi}_{n}(0))
       \big[\Omega^*_n(z_j)-\frac{\kappa_{n+1}}{\kappa_{n}}\Theta^*_n(z_j)\big]
    - \sum_{k \neq j}\frac{\rho_k}{2V(z_k)}\frac{\dot{z}_j-\dot{z}_k}{z_j-z_k}
  \\ \times
      \left[ z_k\Theta^*_n(z_k)
                \big[\Omega^*_n(z_j)-\frac{\kappa_{n+1}}{\kappa_{n}}\Theta^*_n(z_j)\big]
            -z_j\Theta^*_n(z_j)
                \big[\Omega^*_n(z_k)-\frac{\kappa_{n+1}}{\kappa_{n}}\Theta^*_n(z_k)\big] \right]
    \bigg\} .
  \label{ops_Schl:c}
\end{multline}

The fact that the deformation equations satisfy the Schlesinger system of 
partial differential equations should be of no great surprise as the isomonodromic
properties of the regular semi-classical weights are quite transparent. 
\begin{proposition}
The monodromy matrix $ M_j $, $ j=1,\ldots,m,\infty $ defined by the analytic 
continuation of $ Y_n $ around a closed loop enclosing the singularity $ z_j $
\begin{equation}
   \left .Y_n \right|_{z_j + \delta e^{2\pi i}} =
   \left .Y_n \right|_{z_j + \delta}M_j,
\end{equation}
is constant with respect to the deformation variable, $ \dot{M}_j=0 $.
\end{proposition}
\begin{proof}
In the neighbourhood of any isolated finite singularity $ |z-z_j| < \Delta $ the 
Carath\'eodory function is a solution of the inhomogeneous first order differential
equation (\ref{ops_FD}) and can be decomposed as
\begin{equation}
   F(z) = f_j(z) + C_jw(z) .
\label{CF_localExp}
\end{equation}
Here $ f_j(z) $ is the unique, holomorphic solution of the inhomogeneous ODE 
in this neighbourhood whose existence is guaranteed by conditions (3) and (4) of 
the definition (\ref{rSC_defn}) and which we express as
\begin{equation}
   f_j(z) = \sum^{\infty}_{l\geq 0} a_{j,l}(z-z_j)^l .
\end{equation}
The second term of (\ref{CF_localExp}) is a solution of the homogeneous form of
(\ref{ops_FD}) and $ C_j $ a coefficient depending only on $ \{z_k,\rho_k\}^{m}_{k=1} $.
We note that the deformation derivative of the Carath\'eodory function can be 
written as 
\begin{equation}
   \dot{F}(z) = \frac{\dot{w}}{w}F(z)
                +\frac{1}{2}\sum^m_{j=1}\rho_j\frac{\dot{z}_j}{z_j}\frac{z+z_j}{z-z_j}F(z_j)
                 -\frac{1}{2}w_0\sum^m_{j=1}\rho_j\frac{\dot{z}_j}{z_j} ,
\end{equation}
under the assumptions that either $ z_k \notin \TT, k=1,\ldots,m $ or if any 
$ z_k \in \TT $ then $ \Re{\rho_k} > 0 $, to ensure the integral (\ref{ops_Cfun})
is uniformly and absolutely convergent. Differentiating (\ref{CF_localExp}) with 
respect to $ t $ then implies the relation
\begin{multline}
   \dot{C}_j(z-z_j)^{\rho_j}\prod^m_{k\neq j}(z-z_k)^{\rho_k}
  + \rho_j\frac{\dot{z}_j}{z-z_j}\left[ f_j(z)-F(z_j) \right] \\
   + \sum^{\infty}_{l\geq 0}(\dot{a}_{j,l}-(l+1)a_{j,l+1})(z-z_j)^l
   + \sum^m_{l\neq j}\rho_l\frac{\dot{z}_l}{z-z_l}[f_j(z)-F(z_l)]
     +  \sum^m_{l}\rho_l\frac{\dot{z}_l}{2z_l}[w_0-F(z_l)]= 0 .
\end{multline}
Making the assumption $ \Re{\rho_j} > 0 $ we can use the equality $ F(z_j)=f_j(z_j) $
and consequently under condition (3) of (\ref{rSC_defn}) we have
\begin{equation}
  \dot{C}_j(z-z_j)^{\rho_j} + \text{analytic function of $ z $} = 0 ,
\end{equation}
for $ |z-z_j| < \min\{\Delta,|z_j-z_k|,||z_j|-1|\} $. 
This implies $ \dot{C}_j=0 $ by condition (4). Using this fact and that the monodromy 
matrix is given by
\begin{equation}
  M_j = \begin{pmatrix}
              1 & C_j(1-e^{-2\pi i\rho_j}) \cr
              0 & e^{-2\pi i\rho_j} \cr
        \end{pmatrix} ,
\end{equation}
the proof is concluded.
\end{proof}

\begin{remark}
The monodromy matrices are all upper triangular, which is consistent with the fact 
that we are dealing with classical solutions of the Schlesinger systems. Also, they
are independent of $ n $, which is to say they are preserved under the iteration
$ n \to n+1 $. 
\end{remark}

\subsection*{Acknowledgments}
This research has been supported by the Australian Research Council.
Our manuscript has benefited from the critical
reading by Alphonse Magnus and we thank him, Mourad Ismail and Percy Deift for their
advice and suggestions.

\bibliographystyle{amsplain}
\bibliography{moment,nonlinear,random_matrices}

\def\cprime{$'$} \def\cprime{$'$} \def\cydot{\leavevmode\raise.4ex\hbox{.}}
  \def\cprime{$'$} \def\cprime{$'$} \def\cprime{$'$} \def\cprime{$'$}
  \def\cprime{$'$} \def\cprime{$'$} \def\cprime{$'$} \def\cprime{$'$}
  \def\cprime{$'$} \def\cdprime{$''$} \def\cydot{\leavevmode\raise.4ex\hbox{.}}
  \def\cydot{\leavevmode\raise.4ex\hbox{.}} \def\cprime{$'$}
  \def\cydot{\leavevmode\raise.4ex\hbox{.}} \def\cprime{$'$} \def\cprime{$'$}
  \def\cprime{$'$} \def\cprime{$'$} \def\cprime{$'$} \def\cprime{$'$}
  \def\cprime{$'$} \def\cprime{$'$} \def\cprime{$'$} \def\cprime{$'$}
  \def\cprime{$'$} \def\cprime{$'$} \def\cydot{\leavevmode\raise.4ex\hbox{.}}
  \def\cprime{$'$} \def\cprime{$'$} \def\cprime{$'$} \def\cprime{$'$}
  \def\cprime{$'$} \def\cprime{$'$} \def\cprime{$'$} \def\cprime{$'$}
  \def\cprime{$'$} \def\cprime{$'$} \def\cprime{$'$} \def\cprime{$'$}
  \def\cprime{$'$} \def\cprime{$'$} \def\cprime{$'$} \def\cprime{$'$}
  \def\cprime{$'$} \def\cprime{$'$} \def\cprime{$'$} \def\cprime{$'$}
  \def\cprime{$'$} \def\cprime{$'$} \def\cprime{$'$} \def\cprime{$'$}
  \def\cprime{$'$} \def\cprime{$'$} \def\cydot{\leavevmode\raise.4ex\hbox{.}}
  \def\cprime{$'$} \def\cprime{$'$} \def\cprime{$'$} \def\cprime{$'$}
  \def\cprime{$'$} \def\cprime{$'$} \def\cprime{$'$} \def\cprime{$'$}
  \def\cprime{$'$} \def\cprime{$'$} \def\cprime{$'$} \def\cprime{$'$}
  \def\cprime{$'$} \def\cprime{$'$} \def\cprime{$'$}
  \def\cydot{\leavevmode\raise.4ex\hbox{.}} \def\cprime{$'$} \def\cprime{$'$}
  \def\cprime{$'$} \def\cprime{$'$} \def\cprime{$'$} \def\cprime{$'$}
  \def\cprime{$'$} \def\cprime{$'$} \def\cprime{$'$} \def\cprime{$'$}
  \def\cprime{$'$} \def\cprime{$'$} \def\cprime{$'$} \def\cprime{$'$}
  \def\cprime{$'$} \def\cprime{$'$} \def\cprime{$'$} \def\cprime{$'$}
  \def\cprime{$'$} \def\cprime{$'$} \def\cprime{$'$} \def\cprime{$'$}
  \def\cprime{$'$} \def\cprime{$'$} \def\cprime{$'$} \def\cprime{$'$}
  \def\cprime{$'$} \def\cprime{$'$} \def\cprime{$'$} \def\cprime{$'$}
  \def\cprime{$'$} \def\cprime{$'$} \def\cydot{\leavevmode\raise.4ex\hbox{.}}
  \def\cydot{\leavevmode\raise.4ex\hbox{.}}
  \def\cydot{\leavevmode\raise.4ex\hbox{.}}
  \def\cydot{\leavevmode\raise.4ex\hbox{.}}
  \def\cydot{\leavevmode\raise.4ex\hbox{.}}
  \def\cydot{\leavevmode\raise.4ex\hbox{.}}
  \def\cydot{\leavevmode\raise.4ex\hbox{.}}
  \def\cydot{\leavevmode\raise.4ex\hbox{.}} \def\cprime{$'$} \def\cprime{$'$}
\providecommand{\bysame}{\leavevmode\hbox to3em{\hrulefill}\thinspace}
\providecommand{\MR}{\relax\ifhmode\unskip\space\fi MR }
\providecommand{\MRhref}[2]{%
  \href{http://www.ams.org/mathscinet-getitem?mr=#1}{#2}
}
\providecommand{\href}[2]{#2}
\begin{thebibliography}{10}

\bibitem{La_1972b}
\emph{Oeuvres de {L}aguerre. {T}ome {I}}, Chelsea Publishing Co., Bronx, N.Y.,
  1972, Alg\`ebre. Calcul int\'egral, R\'edig\'ees par Ch. Hermite, H.
  Poincar\'e et E. Rouch\'e, R\'eimpression de l'\'edition de 1898. \MR{52
  \#13292}

\bibitem{Ba_1983a}
V.~M. Badkov, \emph{Asymptotic properties of orthogonal polynomials},
  Constructive function theory '81 (Varna, 1981), Publ. House Bulgar. Acad.
  Sci., Sofia, 1983, pp.~21--27. \MR{85c:42026}

\bibitem{Ba_1983b}
\bysame, \emph{Uniform asymptotic representations of orthogonal polynomials},
  Trudy Mat. Inst. Steklov. \textbf{164} (1983), 6--36, Orthogonal series and
  approximations of functions. \MR{86d:42023}

\bibitem{Ba_1992}
\bysame, \emph{Asymptotic and extremal properties of orthogonal polynomials
  with singularities in the weight}, Trudy Mat. Inst. Steklov. \textbf{198}
  (1992), 41--88. \MR{95h:42027}

\bibitem{BDJ_1999}
J.~Baik, P.~Deift, and K.~Johansson, \emph{On the distribution of the length of
  the longest increasing subsequence of random permutations}, J. Amer. Math.
  Soc. \textbf{12} (1999), no.~4, 1119--1178. \MR{2000e:05006}

\bibitem{BDJ_2000a}
\bysame, \emph{On the distribution of the length of the second row of a {Y}oung
  diagram under {P}lancherel measure}, Geom. Funct. Anal. \textbf{10} (2000),
  no.~4, 702--731. \MR{2001m:05258a}

\bibitem{BDMcLMZ_2001}
J.~Baik, P.~Deift, K.T.-R. McLaughlin, P.~Miller, and X.~Zhou, \emph{Optimal
  tail estimates for directed last passage site percolation with geometric
  random variables}, Adv. Theor. Math. Phys. \textbf{5} (2001), no.~6,
  1207--1250. \MR{1 926 668}

\bibitem{BKMcLM_2004}
J.~Baik, T.~Kriecherbauer, K.~T.-R. McLaughlin, , and P.~D. Miller,
  \emph{{Uniform Asymptotics for Polynomials Orthogonal With Respect to a
  General Class of Discrete Weights and Universality Results for Associated
  Ensembles}}.

\bibitem{Ba_2001}
Jinho Baik, \emph{Riemann-{H}ilbert problems for last passage percolation},
  Recent developments in integrable systems and Riemann-Hilbert problems
  (Birmingham, AL, 2000), Contemp. Math., vol. 326, Amer. Math. Soc.,
  Providence, RI, 2003, pp.~1--21. \MR{1 989 002}

\bibitem{Ba_1990}
W.~C. Bauldry, \emph{Estimates of asymmetric {F}reud polynomials on the real
  line}, J. Approx. Theory \textbf{63} (1990), no.~2, 225--237. \MR{92c:33008}

\bibitem{Ba_1960}
Glen Baxter, \emph{Polynomials defined by a difference system}, Bull. Amer.
  Math. Soc. \textbf{66} (1960), 187--190. \MR{22 \#2827}

\bibitem{Ba_1961}
\bysame, \emph{Polynomials defined by a difference system}, J. Math. Anal.
  Appl. \textbf{2} (1961), 223--263. \MR{23 \#A3421}

\bibitem{BR_1994}
S.~Belmehdi and A.~Ronveaux, \emph{Laguerre-{F}reud's equations for the
  recurrence coefficients of semi-classical orthogonal polynomials}, J. Approx.
  Theory \textbf{76} (1994), no.~3, 351--368. \MR{95f:42038}

\bibitem{BEH_2003a}
M.~Bertola, B.~Eynard, and J.~Harnad, \emph{Differential systems for
  biorthogonal polynomials appearing in 2-matrix models and the associated
  {R}iemann-{H}ilbert problem}, Comm. Math. Phys. \textbf{243} (2003), no.~2,
  193--240. \MR{2004i:34237}

\bibitem{BEH_2003b}
\bysame, \emph{Partition functions for matrix models and isomonodromic tau
  functions}, J. Phys. A \textbf{36} (2003), no.~12, 3067--3083, Random matrix
  theory. \MR{2004e:82015}

\bibitem{BC_1990}
S.~S. Bonan and D.~S. Clark, \emph{Estimates of the {H}ermite and the {F}reud
  polynomials}, J. Approx. Theory \textbf{63} (1990), no.~2, 210--224.
  \MR{92c:33007}

\bibitem{CMN_1991}
Giuliana Criscuolo, Giuseppe Mastroianni, and Paul Nevai, \emph{Associated
  generalized {J}acobi functions and polynomials}, J. Math. Anal. Appl.
  \textbf{158} (1991), no.~1, 15--34. \MR{92h:33018}

\bibitem{De_1999}
P.~A. Deift, \emph{Orthogonal polynomials and random matrices: a
  {R}iemann-{H}ilbert approach}, New York University Courant Institute of
  Mathematical Sciences, New York, 1999. \MR{2000g:47048}

\bibitem{EMN_1994}
Tam{\'a}s Erd{\'e}lyi, Alphonse~P. Magnus, and Paul Nevai, \emph{Erratum:
  ``{G}eneralized {J}acobi weights, {C}hristoffel functions, and {J}acobi
  polynomials''}, SIAM J. Math. Anal. \textbf{25} (1994), no.~5, 1461.
  \MR{95f:33011b}

\bibitem{EN_1992}
Tam{\'a}s Erd{\'e}lyi and Paul Nevai, \emph{Generalized {J}acobi weights,
  {C}hristoffel functions, and zeros of orthogonal polynomials}, J. Approx.
  Theory \textbf{69} (1992), no.~2, 111--132. \MR{93c:33004}

\bibitem{FIK_1991}
A.~S. Fokas, A.~R. It{\cydot{s}}, and A.~V. Kitaev, \emph{Discrete {P}ainlev\'e
  equations and their appearance in quantum gravity}, Comm. Math. Phys.
  \textbf{142} (1991), no.~2, 313--344. \MR{93a:58080}

\bibitem{FIK_1992}
\bysame, \emph{The isomonodromy approach to matrix models in $2${D} quantum
  gravity}, Comm. Math. Phys. \textbf{147} (1992), no.~2, 395--430.
  \MR{93h:81115}

\bibitem{rmt_Fo}
P.~J. Forrester, \emph{{L}og {G}ases and {R}andom {M}atrices}, \hfil\break{\tt
  http://www.ms.unimelb.edu.au/$\tilde{\,}$matpjf/matpjf.html}.

\bibitem{FW_2002b}
P.~J. Forrester and N.~S. Witte, \emph{{A}pplication of the $\tau$-function
  theory of {P}ainlev{\'e} equations to random matrices: {PVI}, the {JUE},
  {CyUE}, {cJUE} and scaled limits}, Nagoya Math. J. \textbf{174} (2004),
  29--114.

\bibitem{FW_2004b}
\bysame, \emph{{D}iscrete {P}ainlev\'e equations, {O}rthogonal {P}olynomials on
  the {U}nit {C}ircle and ${N}$-recurrences for averages over ${U(N)}$ --
  {\PVI} $\tau$-functions}, 2004.

\bibitem{GN_1982}
J.~L. Gammel and J.~Nuttall, \emph{Note on generalized {J}acobi polynomials},
  The Riemann problem, complete integrability and arithmetic applications
  (Bures-sur-Yvette/New York, 1979/1980), Lecture Notes in Math., vol. 925,
  Springer, Berlin, 1982, pp.~258--270. \MR{84g:33018}

\bibitem{Ge_1958}
B.~Germansky, \emph{An identity in the theory of the generalized polynomials of
  {J}acobi.}, Proc. Amer. Math. Soc. \textbf{9} (1958), 953--956. \MR{23
  \#A1074}

\bibitem{Ge_1961}
Ya.~L. Geronimus, \emph{Orthogonal polynomials: {E}stimates, asymptotic
  formulas, and series of polynomials orthogonal on the unit circle and on an
  interval}, Authorized translation from the Russian, Consultants Bureau, New
  York, 1961. \MR{24 \#A3469}

\bibitem{Ge_1962}
\bysame, \emph{Polynomials orthogonal on a circle and their applications},
  Translations, {S}er. 1, {V}ol. 3: {S}eries and approximation, vol.~3,
  American Mathematical Society, Providence, R.I., 1962, p.~78.

\bibitem{Ge_1977}
\bysame, \emph{Orthogonal polynomials}, Amer. Math. Soc. Transl. Series 2, vol.
  108, American Mathematical Society, Providence, Rhode Island, 1977.

\bibitem{Hi_1996}
M.~Hisakado, \emph{Unitary matrix models and {P}ainlev\'e {I}{I}{I}}, Mod.
  Phys. Lett. A \textbf{11} (1996), no.~38, 3001--3010. \MR{97h:81200}

\bibitem{IR_1992}
M.~E.~H. Ismail and R.~W. Ruedemann, \emph{Relation between polynomials
  orthogonal on the unit circle with respect to different weights}, J. Approx.
  Theory \textbf{71} (1992), no.~1, 39--60. \MR{94f:42028}

\bibitem{IW_2001}
M.~E.~H. Ismail and N.~S. Witte, \emph{Discriminants and functional equations
  for polynomials orthogonal on the unit circle}, J. Approx. Theory
  \textbf{110} (2001), no.~2, 200--228. \MR{2002e:33011}

\bibitem{IKF_1991}
A.~R. It{\cydot{s}}, A.~V. Kitaev, and A.~S. Fokas, \emph{Matrix models of
  two-dimensional quantum gravity, and isomonodromic solutions of {P}ainlev\'e
  ``discrete equations''}, Zap. Nauchn. Sem. Leningrad. Otdel. Mat. Inst.
  Steklov. (LOMI) \textbf{187} (1991), no.~Differentsialnaya Geom. Gruppy Li i
  Mekh. 12, 3--30, 171, 174. \MR{93e:81114}

\bibitem{JNT_1989}
W.~B. Jones, O.~Nj{\aa}stad, and W.~J. Thron, \emph{Moment theory, orthogonal
  polynomials, quadrature, and continued fractions associated with the unit
  circle}, Bull. London Math. Soc. \textbf{21} (1989), no.~2, 113--152.
  \MR{90e:42027}

\bibitem{Ma_1994}
A.~P. Magnus, \emph{Painlev\'e equations for semi-classical recurrence
  coefficients}, {\tt CA/9409228}, 1994.

\bibitem{Ma_1995a}
\bysame, \emph{Painlev\'e-type differential equations for the recurrence
  coefficients of semi-classical orthogonal polynomials}, Proceedings of the
  Fourth International Symposium on Orthogonal Polynomials and their
  Applications (Evian-Les-Bains, 1992), vol.~57, 1995, pp.~215--237.
  \MR{96f:42027}

\bibitem{Ma_1999}
\bysame, \emph{Freud's equations for orthogonal polynomials as discrete
  {P}ainlev\'e equations}, Symmetries and integrability of difference equations
  (Canterbury, 1996), London Math. Soc. Lecture Note Ser., vol. 255, Cambridge
  Univ. Press, Cambridge, 1999, pp.~228--243. \MR{2000k:42036}

\bibitem{MR_1992}
F.~Marcell{\'a}n and I.~A. Rocha, \emph{On semiclassical linear functionals:
  integral representations}, Proceedings of the Fourth International Symposium
  on Orthogonal Polynomials and their Applications (Evian-Les-Bains, 1992),
  vol.~57, 1995, pp.~239--249. \MR{96f:42028}

\bibitem{MR_1998}
\bysame, \emph{Complex path integral representation for semiclassical linear
  functionals}, J. Approx. Theory \textbf{94} (1998), no.~1, 107--127.
  \MR{99h:42047}

\bibitem{Ma_1987}
P.~Maroni, \emph{Prol\'egom\`enes \`a l'\'etude des polyn\^omes orthogonaux
  semi-classiques}, Ann. Mat. Pura Appl. (4) \textbf{149} (1987), 165--184.
  \MR{89c:33016}

\bibitem{MV_1997}
G.~Mastroianni and P.~V{\'e}rtesi, \emph{Some applications of generalized
  {J}acobi weights}, Acta Math. Hungar. \textbf{77} (1997), no.~4, 323--357.
  \MR{99e:42037}

\bibitem{NEM_1994}
Paul Nevai, Tam{\'a}s Erd{\'e}lyi, and Alphonse~P. Magnus, \emph{Generalized
  {J}acobi weights, {C}hristoffel functions, and {J}acobi polynomials}, SIAM J.
  Math. Anal. \textbf{25} (1994), no.~2, 602--614. \MR{95f:33011a}

\bibitem{Si_2004}
B.~Simon, \emph{Orthogonal polynomials on the {U}nit {C}ircle}, 2004,
  forthcoming.

\bibitem{ops_Sz}
G.~Szeg{\"o}, \emph{Orthogonal polynomials}, third ed., Colloquium Publications
  {\bf 23}, American Mathematical Society, Providence, Rhode Island, 1967.

\bibitem{Ve_1997}
P.~V{\'e}rtesi, \emph{On the zeros of generalized {J}acobi polynomials}, Ann.
  Numer. Math. \textbf{4} (1997), no.~1-4, 561--577, The heritage of P. L.
  Chebyshev: a Festschrift in honor of the 70th birthday of T. J. Rivlin.
  \MR{98c:33013}

\bibitem{Ve_1999b}
\bysame, \emph{Asymptotics of derivatives of orthogonal polynomials based on
  generalized {J}acobi weights. {S}ome new theorems and applications}, New
  developments in approximation theory (Dortmund, 1998), Internat. Ser. Numer.
  Math., vol. 132, Birkh\"auser, Basel, 1999, pp.~329--339. \MR{2001f:42041}

\bibitem{Ve_1999a}
\bysame, \emph{Uniform asymptotics of derivatives of orthogonal polynomials
  based on generalized {J}acobi weights}, Acta Math. Hungar. \textbf{85}
  (1999), no.~1-2, 97--130. \MR{2001b:33017}

\bibitem{Ve_2001}
\bysame, \emph{Orthogonal polynomials based on varying {J}acobi-type weights},
  Studia Sci. Math. Hungar. \textbf{38} (2001), 385--402. \MR{2002k:42055}

\end{thebibliography}

\end{document}